

\input amstex
\documentstyle{amsppt}
\magnification=1200
\catcode`\@=11
\redefine\logo@{}
\catcode`\@=13

\define \bn{\Bbb N}
\define \bz{\Bbb Z}
\define \bq{\Bbb Q}
\define \br{\Bbb R}
\define \bc{\Bbb C}
\define \bh{\Bbb H_2}
\define \bp{\Bbb P}

\define \M{{\Cal M}}
\define\Ha{{\Cal H}}
\define\La{{\Cal L}}
\define\geg{{\goth g}}
\define\0o{{\overline 0}}
\define\1o{{\overline 1}}

\define\gi{\Gamma_{\infty}}

\define\hh{{\goth h}}
\define\cc{{\goth c}}

\define\rr{{\goth r}}
\define\ch{\text{ch\ }}
\define\cha{{\text{ch}_\pm\ }}
\define\io{{\overline i}}
\define\mult{\text{mult}}
\define\re{\text{re}}
\define\im{\text{im}}
\define\pr{\text{pr}}
\TagsOnRight

\document

\topmatter

\title
Siegel automorphic form corrections of some Lorentzian
Kac--Moody Lie algebras
\endtitle

\author
Valeri A. Gritsenko \footnote{Supported by
SFB 170 ``Geometrie und Analysis''. \hfill\hfill} and
Viacheslav V. Nikulin \footnote{Partially supported by
Grant of Russian Fund of Fundamental Research;
Grant of ISF MI6000; Grant of ISF and
Russian Government MI6300; Grant of AMS; SFB 170
``Geometrie und Analysis''. \hfill\hfill}
\endauthor

\address
St. Petersburg Department Steklov Mathematical Institute
Fontanka 27,
\newline
${}\hskip 8pt $
191011 St. Petersburg,  Russia
\endaddress
\email
gritsenk\@lomi.spb.su;\ \ gritsenk\@cfgauss.uni-math.gwdg.de
\endemail

\address
Steklov Mathematical Institute,
ul. Vavilova 42, Moscow 117966, GSP-1, Russia.
\endaddress
\email
slava\@nikulin.mian.su
\endemail

\abstract
We find automorphic form corrections which are generalized
Lorentzian Kac--Moody superalgebras without odd real simple roots
(see R. Borcherds \cite{Bo1} -- \cite{Bo7},
V. Kac \cite{Ka1} -- \cite{Ka3}, R. Moody \cite{Mo} and
\S~6 of this paper) for two elliptic Lorent\-zian
Kac--Moody algebras of the rank $3$ with a lattice Weyl vector,
and calculate multiplicities of their simple and arbitrary imaginary
roots (see an appropriate general setting in \cite{N5}).
These Kac--Moody algebras are defined by hyperbolic
(i.e. with exactly one negative square) symmetrized generalized
Cartan matrices
$$
G_1\ =\
\left(
\matrix\hphantom{-} 2 & -2 & -2 \\
       -2 & \hphantom{-} 2 & -2 \\
       -2 & -2 & \hphantom{-} 2
\endmatrix
\right)\hskip30pt \text{and} \hskip 30pt
G_2\ =\ \left(
\matrix
\hphantom{-}4 &  -4 & -12& -4 \\
 -4 &\hphantom{-} 4 & -4 & -12\\
-12 &-4 & \hphantom{-} 4 & - 4 \\
 -4 & -12 & -4 & \hphantom{-}4
\endmatrix
\right)
$$
of the rank $3$.
Both these algebras have elliptic type (i.e. their Weyl groups have
fundamental polyhedra of finite volume in corresponding
hyperbolic spaces) and have a lattice Weyl vector. The correcting
automoprhic forms are Siegel modular forms. The form corresponding to
$G_1$ is the classical Siegel cusp form $\Delta_5 (Z)$ of weight $5$
which is the product of ten even theta-constants.
In particular we find an infinite product formula for $\Delta_5(Z)$.

To find  the correcting automorphic form for $G_2$,
we use the arithmetic lifting of Jacobi forms
on domains of type IV which was constructed
by first author \cite{G1} and \cite{G2}.
It is proven by second author \cite{N5} that the set of elliptic
Lorentzian Kac--Moody algebras with a lattice Weyl
vector is finite
(but may be extremely big). Conjecturally all of them
have automorphic form corrections.
\endabstract

\rightheadtext
{Siegel automorphic forms and Kac--Moody algebras}

\leftheadtext{V. Gritsenko and  V. Nikulin}

\endtopmatter

\document

\head
\S~0. Introduction
\endhead

First, we explain results of this paper formally.
We find automorphic form corrections which are generalized
Lorentzian Kac--Moody superalgebras without odd real simple roots
(see works of R. Borcherds \cite{Bo1} -- \cite{Bo7},
V. Kac \cite{Ka1} -- \cite{Ka3}, R.~Moody \cite{Mo} and
\S~6 of this paper) for two elliptic Lorent\-zian
Kac--Moody algebras, and
calculate multiplicities of
their simple and arbitrary imaginary roots.
One can find an appropriate general setting on this subject in \cite{N5}.
The first Kac--Moody algebra has the symmetric generalized Cartan matrix
$$
G_1=(\delta_i, \delta_j)=
\left(
\matrix \hphantom{-}2 & -2 & -2 \\
       -2 & \hphantom{-} 2 & -2 \\
       -2 & -2 & \hphantom{-} 2
\endmatrix
\right).
\tag{0.1}
$$
The second one has the symmetrized generalized Cartan matrix
$$
G_2=(\delta_i,\delta_j)=
\left(
\matrix
\hphantom{-}4 &  -4 & -12& -4 \\
 -4 & \hphantom{-}4 & -4 & -12\\
-12 &-4 & \hphantom{-} 4 & - 4 \\
 -4 & -12 & -4 & \hphantom{-}4
\endmatrix
\right).
\tag{0.2}
$$
For both these elliptic Lorentzian Kac--Moody
algebras, the correcting automorphic forms are holomorphic
modular forms on the Siegel upper half-plane.
They are invariant (with a multiplier
system) with respect to $Sp_4(\bz)$ or  the so-called
paramodular group corresponding  to  the polarization $(1,2)$.
The first form  is especially interesting because it is
the classical cusp  form $\Delta_5(Z)$ of weight $5$
which is the product of all even theta-constants (see \thetag{1.2}).
It is well known (see \cite{F}) that
$\Delta_5(Z)^2=F_{10}(Z)$ where $F_{10}(Z)$ is one of the standard
generators for the ring of Siegel modular forms with respect
to $Sp_4(\bz )$.
As a result, we find the corresponding product
formula for this classical Siegel cusp form $\Delta_5(Z)$.

To construct the product formula for $\Delta_5(Z)$,
we use the  Maass's lifting for $Sp_4(\bz)$ (see
\cite{M2}, \cite{EZ}). To find the product formula for correcting
automorphic form of the algebra of $G_2$, we use
the  generalization of this lifting
to the  cases of  the  paramodular subgroups of $Sp_4(\bq)$
and the orthogonal  groups of signature $(2,n)$
(see  \cite{G1} -- \cite{G3}).

We should say that  first examples of automorphic
form corrections
of Lorent\-zian Kac--Moody algebras
were found by R. Borcherds (see \cite{Bo1}---\cite{Bo7}). Mainly,
his examples are connected with
automoprhic forms on the complex domains of large  dimension
(for example the dimension  $26$)
and the arithmetic of the Leech lattice.

Our examples above are first examples of
automorphic form corrections of elliptic (or parabolic)
Lorentzian Kac--Moody algebras
(we use terminology from \cite{N5})
when correcting automorphic forms are Siegel modular forms
(necessarily of the genus $2$).
Moreover, we show that there are
very classical automoprhic forms like
$\Delta_5(Z)$ above which give the automorphic
form corrections of some Lorentzian Kac--Moody algebras.
In forthcoming publication we hope to give  other  examples
of automorphic forms which correct the corresponding Lorentzian
Kac--Moody algebras.

\vskip5pt

Now we want to present our results in some details.
We restrict considering  the case \thetag{0.1}.

Hyperbolic symmetric matrix $G_1$ defines an integral hyperbolic
symmetric bilinear form (i.e. hyperbolic lattice)
on the $\bz$-module $M_{II}=\bz \delta_1\oplus \bz \delta_2\oplus
\bz \delta_3$. Reflections $s_{\delta_i}$ in $\delta_i \in M_{II}$
generate a reflection subgroup of finite index $W\subset O(M_{II})$
which is discrete in corresponding hyperbolic space $V^+(M_{II})/\br_{++}$
and has a fundamental polyhedron $\M_{II}$ of finite volume
which is bounded by half-spaces
orthogonal to $\delta_i$. Here $V^+(M_{II})$ is a half-cone of the cone
$V(M_{II})=\{x \in M_{II}\otimes \br \ |\ (x,x)<0\}$. The half-cone
$V^+(M_{II})$ is choosen uniquely by the following condition which is
equivalent to finiteness of volume of $\M_{II}$:
$$
V^+(M_{II})\subset \br_+\delta_1+\br_+\delta_2+\br_+\delta_3.
$$
Moreover, $\br_+\M_{II}$ is the dual cone
$$
\br_+\M_{II} =\{x \in M_{II}\otimes \br\ |\ (x, \delta_i)\le 0,\ i=1,2,3\}
\subset \overline{V^+(M_{II})}
$$
to  $\br_+\delta_1+\br_+\delta_2+\br_+\delta_3$.
We denote $P(\M_{II})=\{\delta_1, \delta_2, \delta_3\}\subset M_{II}$.
The group $W$
and the set $P(\M_{II})$ of vectors orthogonal to the fundamental polyhedron
$\M_{II}$ of $W$ have a {\it lattice Weyl vector} which is an element
$\rho \in M_{II}\otimes \bq$ with the property
$$
(\rho , \delta_i)=-(\delta_i, \delta_i)/2.
$$
Evidently, $\rho =(\delta_1+\delta_2+\delta_3)/2$. It is shown in \cite{N5}
(in fact, it easily follows from results of \cite{N1} and \cite{N2}) that
the set of symmetrized generalized Cartan matrices with these properties
(i. e. finiteness of volume of the fundamental polyhedron of $W$ and
existence of a lattice Weyl vector for the set of vectors
$P(\M)$ orthogonal to a fundamental polyhedron $\M$ of $W$)
is finite (up to multiplication by a constant) for rank $\ge 3$.
 Thus, we are considering one of the finite set of cases.

The cone $V^+(M_{II})$ defines the corresponding complexified cone
$\Omega (V^+(M_{II}))=M_{II}\otimes \br + iV^+(M_{II})$.
We consider
extended lattice $L$ which is an orthogonal sum of $M_{II}$ and
unimodular even hyperbolic plane $U$. One can consider
$\Omega (V^+(M_{II}))$ as a cusp of the domain
${\Cal H}_+$ of type IV which is one of two connected components of
$$
{\Cal H} = \{ \bc \omega \subset L\otimes \bc\ |\  (\omega, \omega)=0,
\ (\omega, \overline{\omega})<0\}.
$$
Using natural identification of Siegel half-space
${\Bbb H}_2$ and ${\Cal H}_+$ and natural isomorphism
$Sp_4(\bz)/\{\pm E_4\}\subset SO^+_\br(L)$,
we can consider $\Delta_5(Z)$
as an automophic form with respect to $SO^+_\br(L)$.
We show that
$\frac{1}{64}
\Delta_5(Z)$ is anti-invariant
with respect to $W$, has integral Fourier coefficients and its
Fourier coefficient corresponding to $\rho$ is equal to $1$.
More exactly, we show (see Theorem 2.3) that
$$
\hskip-250pt {\tsize\frac{1}{64}} \Delta_5(z)=
$$
$$
\sum_{w\in W}{ \det(w)
\left(\exp(-\pi i (w(\rho),z))
- \hskip-10pt  \sum_{ a \in M_{II} \cap \br_{++}\M_{II}}
{ \hskip-10pt   m(a)\exp(-\pi i (w(\rho+a),z))}\right)},
\tag{0.3}
$$
where $z=z_3f_2+z_2f_3+z_1f_{-2} \in M_{II}
\otimes \br +i V^+(M_{II})$ and
$m(a) \in \bz$;
for any primitive $a_0 \in M_{II}\cap \br_{++}\M_{II}$
with $(a_0,a_0)=0$, we have the identity
$$
1- \sum_{t \in \bn}{m(ta_0)q^t}=
\prod_{k \in \bn}{(1-q^k)^{\tau(ka_0)}}
$$
of power series of one variable $q$, where all $\tau(ka_0)=9$.

The matrix \thetag{0.1} is the Gram matrix of elements $P(\M_{II})$.
We use Fourier coefficients in \thetag{0.3} to
extend the set $P(\M)$ to the set
${}_s\Delta = {}_s\Delta^{\re} \cup {}_s\Delta^{\im}$ where
${}_s\Delta^{\re}=P(\M )$ and
$$
\split
{}_s\Delta^{\im}&=
\{m(a)a\ |\ a\in M_{II}\cap \br_{++}\M_{II},\ (a,a)<0 \}\\
&\cup\{\tau(a)a\ |\ a\in M_{II}\cap \br_{++}\M_{II},\ (a,a)=0 \}.
\endsplit
$$
One can use Gram matrix of ${}_s\Delta$ to construct a generalized
Kac--Moody superalgebra  $\geg$ with the set of simple roots
${}_s\Delta$. This construction is similar to  \cite{Bo1}
(see also \cite{Bo2}),
 \cite{Ka1}---\cite{Ka3} and \cite{Mo}.
See \S~3 and \S~6   for details.
The Weyl--Kac denominator function $\Phi (z)$ of this
Kac--Moody superalgebra
is equal to $\frac{1}{64}\Delta_5(2z)$ since our construction of $\geg$.
The Kac--Moody superalgebra
$\geg$ is called ``corrected'' using the automoprhic form
$\frac{1}{64}\Delta_5(z)$. The Weyl--Kac denominator function $\Phi (z)$
for $\geg$ has the product formula
$$
\hskip-10pt\Phi (z)= \sum_{w\in W }{\det(w)
\hskip-2pt \left(\hskip-2pt \exp(-2\pi i (w(\rho),z))
- \hskip-15pt \sum_{ \hskip-10pt a \in M_{II} \cap \br_{++}\M_{II}}
{\hskip-15pt m(a)\exp(-2\pi i (w(\rho+a),z))}\hskip-2pt \right)}\hskip-2pt =
$$
$$
=\exp({-2\pi i(\rho,z)})\prod_{\alpha \in \Delta_+}
{\left( 1-\exp({-2\pi i (\alpha , z)})\right)^
{\mult~\alpha}}.
\tag{0.4}
$$
Thus, we should have similar product formula for $\Delta_5(z)$. Using
automoprhic properties of $\Delta_5(z)$ and Maass's lifting, we
find this formula. In particular, we calculate multiplicities
$\mult~\alpha$. See Theorem 4.1.

\vskip5pt

In \cite{G1} -- \cite{G3}, the arithmetic
lifting of Jacobi forms on domains of IV type is constructed.
This is a  generalization of
Maass's lifting \cite{M2}. We use this lifting to get similar results
for the Kac--Moody algebra corresponding to the matrix (0.2). See \S~5.

\vskip5pt

Both our examples (0.1) and (0.2) belong to one series of
elliptic Lorentzian Kac--Moody algebras related with
reflection groups on hyperbolic plane which were classified in
\cite{N3}. They are groups $W^{(2)}(S)$ generated by reflections in
all elements with square $2$ of hyperbolic (i.e. of the signature
$(n,1)$) integral lattices $S$ of the rank $3$ such
that index $[O(S):W^{(2)}(S)]$ is finite. We hope to consider
other and may be all examples from this series in further
publications.

\vskip5pt

{\bf Acknowledgment.} This paper was written during our stay
at SFB 170  ``Geometrie und Analysis'' of Mathematical
Institute of Georg-August-University
in  January -- March 1995 and Steklov Mathematical Institute.
We are grateful to these Institutes for hospitality.

\head
\S~1. Siegel modular forms
\endhead

We remind that Siegel modular group $Sp_4(\bz)$ is the group of
of all integral $4\times 4$-
matrices
$$
\qquad M=\pmatrix A&B\\C&D\endpmatrix\
$$
such that $^tM J_4M=J_4$ where $A,B,C,D$ are $2\times 2$-integral matrices
and
$$
\qquad {J_4}=\pmatrix 0 & E_2 \\
                      -E_2 &0\\
\endpmatrix .
$$
Siegel modular form $F$ of the weight $k$ is a holomorphic function
on the Siegel domain
$$
\bh=\{Z={}^tZ\in M_2(\bc),\ Z=X+iY,\ Y>0\}
$$
with the invariance property
$$
F\bigl((AZ+B)(CZ+D)^{-1})=\hbox{det\,}(CZ+D)^k F(Z),\
\qquad M=\pmatrix A&B\\C&D\endpmatrix\in Sp_4(\bz).
$$
These forms define the ring
$$
\frak M (Sp_4(\bz))=\bigoplus_k \frak M_k (Sp_4(\bz))
$$
where $\frak{M}_k (Sp_4(\bz))$ denote the subspace of functions
of the weight $k$.
This ring has four generators
$$
\frak M (Sp_4(\bz))=\bc\,[E_4,\,E_6,\,F_{10},\, F_{12}],
\tag1.1
$$
where  $E_4(Z)$, $E_6(Z)$
are the Eisenstein series and $F_{10}(Z)$ and $F_{12}(Z)$
are cusp forms  (here below index
denote weight). Moreover, $F_{10}(Z)$ is the square of a cusp form
$\Delta_5(Z)$ of weight $5$ with a multiplier system (see \cite{F}).
The function $\Delta_5(Z)$ is given by the product of all even
theta-constants
$$
\Delta_5 (Z)=\prod_{(a,b)}\vartheta_{a,b}(Z),
\tag1.2
$$
where
$$
\vartheta_{a,b}(Z)=\sum_{l\in \bz^2}
\exp{\bigr(\pi i (Z[l+\frac 1{2}a]+{}^tbl)\bigl)}
\qquad\qquad (Z[l]={}^tlZl)
$$
and the product is taken over all vectors $a,b\in (\bz/2\bz)^2$
such that ${}^t ab\equiv 0\mod 2$.
(There are exactly ten different $(a,b)$.)

\vskip5pt

In this paper, it is more appropriate for our purpose to consider
a Siegel modular form as a form with respect to the orthogonal
group $SO(3,2)$.

First, we construct some isomorphisms between symplectic and orthogonal
gro\-ups.

We denote $L_4=\Bbb Z e_1\oplus \Bbb Z e_2 \oplus \Bbb Z e_3
\oplus \Bbb Z e_4$. Any $\bz$-linear map $g:\,L_4\to L_4$ induces a linear map
$\wedge^2 g: L_4\wedge L_4\to L_4\wedge L_4$.

The scalar product (the Pfaffian) on $L_4\wedge L_4$ is defined by
$u\wedge v=(u,v)e_1\wedge e_2\wedge e_3\wedge e_4\in \wedge^4 L$,
$x,y \in L_4\wedge L_4$ . This is
an even unimodular integral symmetric bilinear
form of the signature $(3,3)$ on $L_4\wedge L_4$. It
is invariant with respect to the action of $SL(L_4)$ which was defined
above.

We fix a skew-symmetric form $J_4$ on $L_4$ by the property:
$$
J_4(x,y)e_1\wedge e_2\wedge e_3\wedge e_4=
-x\wedge y\wedge (e_1\wedge e_3+ e_2\wedge e_4).
$$
(Similarly, all elements of $L_4\wedge L_4$ are identified with integral
skew-symmetric bilinear forms on $L_4$.) Therefore
$$
Sp_4(\bz) \cong
\{g:\ L_4\to L_4\ |\  (g\wedge g)(e_1\wedge e_3+ e_2\wedge e_4)
=e_1\wedge e_3+e_2\wedge e_4\}.
$$
Evidently, then $Sp_4(\bz)$ keeps the lattice
$$
L=(e_1\wedge e_3+e_2\wedge e_4)^{\perp}
\subset L_4\wedge L_4,\ \ L\cong U\oplus U\oplus <2>,
$$
where $U$ is an unimodular integral hyperbolic plane, i.e.
a lattice with the quadratic form
$$
\pmatrix 0&-1\\-1&0\endpmatrix.
$$
The $<2>$ is a one dimensional $\bz$-lattice
with the matrix $(2)$.
We fix  a  basis in $L$
$$
(f_1=e_1\wedge e_2,\, f_2=e_2\wedge e_3,\, f_3=e_1\wedge e_3
-e_2\wedge e_4,\,
f_{-2}=e_4\wedge e_1,\, f_{-1}=e_4\wedge e_3).
$$

The real orthogonal group $O_{\br}(L)=O(L\otimes \br)$ acts on a domain
$$
\Cal H^{IV}=\{Z\in \bp(L\otimes \bc)\ |\ (Z,Z)=0,\ (Z,\overline Z)<0\}=
\Cal H_+\cup \overline {\Cal H}_+,
$$
where (using the basis above)
$$
\Cal H_+ =
\{Z=^{t} \bigl((z_2^2-z_1z_3),\,z_3,\,z_2,\,z_1,\,1\bigr)\cdot z_0\in
\Cal H^{IV}\ |\ \hbox{Im}\,(z_1)>0\},
$$
which is the classical homogeneous domain of type IV.
The condition $(Z,\overline Z)<0$ is equivalent to
$$
y_1y_3-y_2^{\ 2} > 0,\qquad\text{where }\ y_i=\hbox{Im}(z_i).
$$
The domain $\Cal H_+$  coincides with the Siegel upper half-plane
$\bh$ if we correspond to the point of $\Cal H_+$ with
the parameters $z_1,z_2,z_3$ above a symmetric matrix
$$
\qquad Z=\pmatrix z_1 &z_2\\z_2&z_3\endpmatrix \in \bh.
$$

The real orthogonal group $O_\br(L)$ has four connected components.
By $O_\br^+(L)$ we denote the subgroup of $O_\br(L)$ of index $2$
consisting  of elements, which leave $\Cal H_+$ invariant. The
kernel of the action of $O_\br^+(L)$ on $\Cal H_+$ is equal to
$\pm E_5$. For our case (odd dimension of $L$), the group
$O_\br^+(L)=\pm E_5SO_\br^+(L)$, where $SO_\br^+(L)$
is the subgroup of elements with real spin-norm equals to one.
We put $O^+(L)=O^+_\br (L) \cap O(L)$ and
$SO^+(L)=O(L)\cap SO_\br^+(L)$.

One can easily find the images of generators of $Sp_4(\bz)$:
$$
\align
\wedge^2\bigl(
\pmatrix 0&E_2\\-E_2&0 \endpmatrix
\bigr)    &=
\pmatrix 0&0&0&0&-1\\0&0&0&-1&0\\0&0&1&0&0\\0&-1&0&0&0\\
-1&0&0&0&0\endpmatrix,
\\
\vspace{3\jot}
\wedge^2\bigl(
\pmatrix E_2&B\\0&E_2\endpmatrix
\bigr)&=
\pmatrix 1&-b_1&2b_2&-b_3&b^2-b_1b_2\\0&1&0&0&b_3
\\0&0&1&0&b_2\\0&0&0&0&b_1\\
0&0&0&0&1\endpmatrix, \qquad\quad
B=\pmatrix b_1&b_2\\b_2&b_3\endpmatrix.
\endalign
$$
We also remark that
$$
\wedge^2\bigl(
\pmatrix U^{*}&0\\0&U\endpmatrix
\bigr)=
 \hbox{det\,}(U)\,\pmatrix
1&0&0&0&0\\0&a^2&-2ab&b^2&0\\0&-ac&ad+bc&-bd&0\\0&c^2&-2cd&d^2&0\\
0&0&0&0&1\endpmatrix,
$$
where
$U=\pmatrix a&b\\c&d\endpmatrix$ and $U^*={}^tU^{-1}$.

Using these considerations, it is not difficult to prove
\proclaim{Lemma 1.1}
The correspondence $\wedge^2$ defines the isomorphism
$$
\wedge^2:\ Sp_4(\bz)/\{\pm E_4\} \to SO^+(L)\cong O^+(L)/\{\pm E_5\}
$$
which gives the commutative diagram
$$
\CD
\bh@>g>>\bh\\
@V VV @V VV\\
\Cal H_+@>g\wedge g>>\Cal H_+
\endCD
$$
for the isomorphism $\bh\to {\Cal H}_+$ above and any $g \in Sp_4(\bz)$.
\endproclaim

\vskip10pt

The modular form $\Delta_5(Z)$ has a non-trivial multiplier system.
An exact formula  for it  was found by H. Maass in \cite{M1}:
$$
\Delta_5(M<Z>)=v(M)\,\hbox{det\,}(CZ+D)^{5}\,\Delta_5(Z),\qquad
M=\pmatrix A&B\\C&D\endpmatrix\in Sp_4(\bz),
\tag1.3
$$
where
$$
\gather
v\bigl(\pmatrix 0&E_2\\-E_2&0 \endpmatrix\bigr)=1, \qquad
v\bigl(\pmatrix E_2&B\\0&E_2\endpmatrix\bigr)=(-1)^{b_1+b_2+b_3},
\tag1.4\\
\vspace{3\jot}
v\bigl(
\pmatrix U^{*}&0\\0&U\endpmatrix
\bigr)=(-1)^{(1+a+d)(1+b+c)+ad}
\tag1.5
\endgather
$$
with $U=\pmatrix a&b\\c&d\endpmatrix\in GL_2(\bz)$,  $U^*={}^tU^{-1}$
and $B=\pmatrix b_1&b_2\\b_2&b_3\endpmatrix\in M_2(\bz)$.
According to \thetag{1.2} and  \thetag{1.4},  $\Delta_5(Z)$ has
Fourier expansion
$$
\Delta_5(Z)=
\sum\Sb n,\,l,\,m\equiv 1\, mod \,2\\4nm-l^2> 0\\
n,m>0\endSb f(n,l,m)
\exp{(\pi i\,(nz_1+lz_2+mz_3))}
\tag1.6
$$
with integral $f(n,l,m)$. (We recall that $\Delta_5(Z)$ is a cusp form.)
Moreover, it is not difficult to
calculate that the Fourier coefficients
$$
f(1,1,1)=64 \ \text{ and }\ \ 64\,|\,f(n,l,m)\ \text{ for all }\ \ n,l,m.
\tag1.7
$$
To see the last property, one should remark that all Fourier coefficients
of the theta-constants  $\vartheta_{a,b}(Z)$
are multiple by $2$ for $a\not=0$.
In \S 4 (see considerations after \thetag{4.8}) we prove  the following
identity for the formal power series
$$
1+\frac{1}{64}\sum_{t \in \bn}{f(1+2t,1,1)q^t}=
\prod_{k \in \bn}{(1-q^k)^9}.
\tag1.8
$$
In what follows we use the properties \thetag{1.3}--\thetag{1.8}
of $\Delta_5(Z)$ for construction of a hyperbolic Kac-Moody Lie algebra.

\head
\S~2.  The hyperbolic lattice $M_{II}$, reflection group $W^{(2)}(M_{II})$
and Fourier expansion of $\Delta_5(Z)$
\endhead

We fix a primitive hyperbolic sublattice $M_0=\bz f_2\oplus \bz f_3 \oplus
\bz f_{-2} \cong U\oplus <2>$ of $L$.
We exend an automorphism $\phi \in O(M_0)$ to be
identical on $(M_0)^\perp_L$. This gives an embedding
$O(M_0)\subset O(L)$. We are interesting in automorphic properties
of $\Delta_5(Z)$ with respect to $O(M_0)$. We remind that
every primitive element $\delta \in M_0$ (or any other integral
symmetric bilinear form) with $(\delta, \delta)>0$ and
$(\delta,\delta)|2(M_0, \delta)$ defines a reflection
$$
s_\delta: x \to x-(2(x,\delta)/(\delta, \delta))\delta, \ \ x \in M_0.
$$
Here, $s_\delta (\delta )=-\delta$ and $s_\delta | \delta^\perp$ is
identical. In particular, every $\delta\in M_0$ with
$(\delta, \delta)=2$ defines the reflection
$$
s_\delta: x \to x-(x,\delta)\delta, \ \ x \in M_0.
$$

\proclaim{Proposition 2.1} We have:
$$
\Delta_5(s_{\delta}(Z))=-\Delta_5(Z)
$$
for reflections $s_\delta$ with respect to elements with square $2$
$$
\delta =\delta_1 =2f_2-f_3,\  \delta_2=2f_{-2}-f_3,
\ \delta_3=f_3 \in M_0
$$
and
$$
\Delta_5(s_\delta(Z))=\Delta_5(Z)
$$
for reflections $s_\delta$ with respect to elements
$$
\delta=f_{-2}-f_2,\  f_2-f_3,\  f_2+f_3 \in M_0
$$
also with square $2$.
\endproclaim

\demo{Proof} To simplify notation, we denote
$$
\widetilde{U}=\wedge^2\bigl(
\pmatrix U^{*}&0\\0&U\endpmatrix
\bigr),\ \ U\in GL(2,\bz).
$$
It is easy to see that
$s_{f_{-2}-f_2}=-\widetilde
{\left(\smallmatrix 0&1\\1&0\endsmallmatrix\right)}$,
$s_{f_3}=-\widetilde{\left(\smallmatrix 1&0\\0&-1
\endsmallmatrix\right)}$,
$s_{f_2-f_3}=-\widetilde{\left(\smallmatrix
1 &-1\\
0 &-1
\endsmallmatrix\right)}$,
therefore
by Maass formula \thetag{1.3}, we get the following    relations
$$
\Delta_5(s_{f_{-2}-f_2}(Z))=\Delta_5(Z),\quad
\Delta_5(s_{f_3}(Z))=-\Delta_5(Z),\quad
\Delta_5(s_{f_3}(Z))=-\Delta_5(Z).
$$
We have
$2f_2-f_3=s_{f_2-f_3}(f_3)$; $2f_{-2}-f_3=s_{f_{-2}-f_2}(2f_2-f_3)$;
$f_2+f_3=s_{f_3}(f_2-f_3)$. From above, it follows (by group's arguments)
that
$\Delta_5(s_{2f_2-f_3}(Z))=-\Delta_5(Z)$,
$\Delta_5(s_{2f_{-2}-f_3}(Z))=-\Delta_5(Z)$,
$\Delta_5(s_{f_2+f_3}(Z))=\Delta_5(Z)$

It follows the statement.
\enddemo

Now we want to interpret Proposition 2.1 more invariantly, using
the automorphism group $O(M_0)$. The lattice $M_0$ is
hyperbolic, i.e. it has signature $(2,1)$. Thus, $M_0$ defines a
cone
$$
V(M_0)=\{x \in M_0\otimes \br\ | \ (x,x)<0\}
$$
which is the union of its two half-cones.
(All general definitions and notations below hold
for an arbitrary hyperbolic (i.e. with signature $(n,1)$)
lattice, and we shall use them later.)
We choose one of this half-cones
uniquely by the condition that
$$
\Omega (V^+(M_0)):=M_0\otimes \br +i V^+(M_0)\subset {\Cal H}_+.
$$
This means that if $z=z_3f_2+z_2f_3+z_1f_{-2} \in \Omega(V^+(M_0))$,
then the point
$Z={}^t((z_2^2-z_1z_3),z_3, z_2, z_1, 1)\cdot z_0 \in \Ha_+$
in notations above.
We denote by $O^+(M_0)$ the subgroup of $O(M_0)$ of index two which fixes
the half-cone $V^+(M_0)$. It is well-known that
the group $O^+(M_0)$ is discrete in the corresponding hyperbolic
space $\La^+(M_0)=V^+(M_0)/\br_{++}$ and has a fundamental domain of finite
volume (below index $++$ denote positive numbers, $+$ denote
non-negative numbers respectively).
Any reflection $s_\delta \in O(M_0)$ with respect to an
element $\delta \in M_0$ with $(\delta, \delta)>0$ is a
reflection in the hyperplane
$$
{\Cal H}_\delta =\{ \br_{++}x \in \La^+(M_0)\ |\ (x, \delta)=0\}
$$
of $\La^+(M_0)$. This maps the half-space
$$
 \Ha^+_\delta =\{ \br_{++}x \in \La^+(M_0)\ |\ (x, \delta)\le 0 \}
$$
to the opposite half-space ${\Cal H}_{-\delta}$ which are both bounded by
the  hyperplane ${\Cal H}_\delta$. Here $\delta \in M_0$ is
called orthogonal to $\Ha_\delta$ and $\Ha^+_\delta$.
All reflections of $M_0$ generate the
reflection subgroup $W(M_0)\subset O^+(M_0)$.

\smallpagebreak

The hyperbolic lattice $M_0$ of Proposition 2.1
is very special, and its automorphism group is
well-known. Below we give results about this group which we need.
These results one certainly can find in \cite{N3}
where the reflection group of the lattice
$M_0$ and some its natural sublattices which are invariant with respect
to its  reflection subgroups of finite index were classified.

First, $O^+(M_0)=W^{(2)}(M_0)$ where  index $2$ denote the subgroup
generated by reflections in all elements of $M_0$ with square $2$.
Similarly, we
denote as $\Delta^{(t)}(K)$
the set of all primitive elements $\delta \in K$ of a lattice $K$
with $(\delta ,\delta)=t$ which define reflections $s_{\delta}$ of $K$,
and as $W^{(t)}(K)$ the group generated by reflections in all
these elements. Thus, $O^+(M_0)$ is generated by
reflections in $\Delta^{(2)}(M_0)$.
An element $\delta \in \Delta^{(2)}(M_0)$ and the corresponding
reflection $s_\delta$ has one of two types:

\vskip5pt

Type ${}\hphantom{I}$I: $(\delta, M_0)= ${}\hphantom{2}$\bz$.

Type II: $(\delta, M_0)=2\bz$.

\vskip5pt

We introduce sublattices $M_I$ and $M_{II}$ which are generated by
all elements $\delta \in \Delta^{(2)} (M_0)$ of the type $I$ and $II$
respectively.
We have
$$
M_I=\{ mf_2+lf_3+nf_{-2} \in M_0 \ | \ m+l+n \equiv 0\mod 2 \}
$$
and
$$
M_{II}=\{mf_2+lf_3+nf_{-2} \in M_0\ |\ m \equiv n \equiv 0\mod 2\}.
$$
An element $\delta\in \Delta^{(2)} (M_0)$ has the type I
(respectively II) if and only if
$\delta \in M_I$ (respectively $\delta \in M_{II}$). It follows that
the subgroup of $O^+(M_0)=W^{(2)}(M_0)$ generated by all reflections
$s_\delta$ of the type I (respectively II) is equal to $W^{(2)}(M_I)$
(respectively $W^{(2)}(M_{II})$).
Obviously, sublattices $M_0, M_I, M_{II}$ are $W^{(2)}(M_0)$-invariant and
both subgroups  $W^{(2)}(M_{I})$ and $W^{(2)}(M_{II})$ are normal in
$W^{(2)}(M_0)$.
The index $[W^{(2)}(M_0):W^{(2)}(M_I)]=2$ and
$[W^{(2)}(M_0):W^{(2)}(M_{II})]=6$. Fundamental polyhedra
$\M=\M_0$, $\M_I$, and $\M_{II}$ for groups
$W^{(2)}(M_0)$, $W^{(2)}(M_I)$ and $W^{(2)}(M_{II})$
respectively are equal to
$\M=\bigcap_{\delta \in P(\M)_{\pr}}\Ha_\delta^+$ where $P(\M)_{\pr}$ are
minimal sets of
primitive elements of $M_0$ with positive square
(with square $2$ for our case). They are
called sets of (primitive) orthogonal vectors to polyhedra and are
equal to
$$
P(\M_0)_{\pr}=\{f_2-f_3, f_{-2}-f_2, f_3\};
$$
$$
P(\M_I)_{\pr}=\{f_2-f_3, f_{-2}-f_2, f_2+f_3\};
$$
and
$$
P(\M_{II})_{\pr}=\{ \delta_1, \delta_2, \delta_3\}.
$$
It follows that the groups $W^{(2)}(M_0)$, $W^{(2)}(M_I)$,
$W^{(2)}(M_{II})$ are generated by reflections
in $P(\M_0)_{\pr}$, $P(\M_I)_{\pr}$ and $P(\M_{II})_{\pr}$ respectively.
We denote
$$
A(P(\M))=\{g \in O^+(M_0)\ |\ g(P(\M))=P(\M) \}
$$
the ``group of symmetries" of a fundamental polyhedron $\M$ and
its set $P(\M)$ of orthogonal vectors.
The group $A(P(\M_0)_{\pr})$ is trivial, $A(P(\M_I)_{\pr})$
has order two and is
generated by $s_{f_3}$, the group $A(P(\M_{II})_{\pr})$ is the
group of symmetries of the right triangle (i.e. it is $S_3$). It
is generated by $s_{f_2-f_3}, s_{f_{-2}-f_2}$.
We can write down the $O^+(M_0)$ as a semi-direct products:
$$
O^+(M_0)=W^{(2)}(M_0)=W^{(2)}(M_I)\rtimes A(P(\M_I)_{\pr})=
W^{(2)}(M_{II})\rtimes A(P(\M_{II})_{\pr}).
$$

Using these considerations,
we can reformulate the invariance properties of
\newline
$\Delta_5(Z)$ of
Proposition 2.1
with respect to the group $O^+(M_0)$ as follows:

\proclaim{Proposition 2.2} We have the following invariance properties
of $\Delta_5(Z)$ with respect to the group
$$
O^+(M_0)=W^{(2)}(M_0)=W^{(2)}(M_I)\rtimes A(P(\M_I)_{\pr})=
W^{(2)}(M_{II})\rtimes A(P(\M_{II})_{\pr}):
$$

{\rm (a)} For a reflection $s_\delta \in O^+(M_0)$ with $(\delta,\delta)=2$
we have
$$
\Delta_5(s_\delta(Z))=
\cases
\Delta_5(s_\delta(Z)),
&\text{if $\delta$ has the type I (i.e $\delta \in M_I$)}, \\
-\Delta_5(s_\delta(Z)),
&\text{if $\delta$ has the type II (i.e $\delta \in M_{II}$)}.
\endcases
$$

{\rm (b)}
$\qquad
\Delta_5(w.a(Z))=\det(a) \Delta_5(Z)\ \text{for\ }
w \in W^{(2)}(M_I),\  a\in A(P(\M_I)_{\pr});
$
\newline
in particular,
$$
\Delta_5(g(Z))=\Delta_5(Z)\ \text{for\ } g\in O^+(M_0)
$$
if and only if $g \in W^{(2)}(M_I)$.

{\rm (c)}
$\qquad
\Delta_5(w.a(Z))=\det (w)\Delta_5(Z)\ \text{for\ }
w \in W^{(2)}(M_{II}), a\in A(P(\M_{II})_{\pr}).
$
\endproclaim

Because of the property (c), the reflection group $W^{(2)}(M_{II})$
will be important for us. We fix the fundamental polyhedron $\M_{II}$
above with the set of orthogonal vectors  $P(\M_{II})=P(\M_{II})_{\pr}=
\{\delta_1, \delta_2, \delta_3\}$.
These elements have the Gram matrix
$$
(\delta_i, \delta_j)=
\left(
\matrix \hphantom{-}2 & -2 & -2 \\
       -2 & \hphantom{-} 2 & -2 \\
       -2 & -2 & \hphantom{-} 2
\endmatrix
\right) .
\tag2.1
$$
Let us consider the cone
$\br_+Q_+=\br_+\delta_1+\br_+\delta_2+\br_+\delta_3$ and the
corresponding dual cone
$\br_+Q_+^\ast=\{x \in M_0\otimes \br\ |\ (x, \delta_i)\le 0\}$.
We know that $\M_{II}\subset \La^+(M_0)$ has finite volume. Since
the cone $\overline{V^+(M_0)}=\overline{V^+(M_0)^\ast}$
is self-dual, this is equivalent
to the sequence of embeddings of cones
$$
\br_+ Q_+^\ast \subset \overline{V^+(M_0)} \subset \br_+Q_+.
\tag2.2
$$
This property (equivalent to finiteness of
volume of $\M_{II}$) is very important. This gives a hope that all
our further considerations may be successful (see \cite{N5}). We will
use this property many times further.

Another very important property of the group $W^{(2)}(M_{II})$ and the
set $P(\M_{II})$ is that they have a {\it lattice Weyl vector} $\rho$.
(See general results on this subject and an explanation why this is
important in \cite{N5}.)
This is the element
$\rho \in M_{II}^\ast$ with the property
$$
(\rho, \delta_i)=-(\delta_i, \delta_i)/2=-1\ \text{for any
$\delta_i\in P(\M_{II})$}.
\tag2.3
$$
By \thetag{2.1},
$$
\rho =(1/2)\delta_1+(1/2)\delta_2+(1/2)\delta_3=f_2-(1/2)f_3+f_{-2}.
\tag2.4
$$
Evidently, $\rho \in \br_+Q_+^\ast$, and by \thetag{2.2},
$\rho \in V^+(M_0)=V^+(M_{II})$. Further we will work with
the lattice $M_{II}$
identifying $M_0\otimes \bq=M_{II}\otimes \bq$.
Using the group $W^{(2)}(M_{II})$, we can rewrite $\Delta_5(Z)$
as follows

\proclaim{Theorem 2.3} {\rm (a)}
$\frac{1}{64} \Delta_5(Z)=$
$$
\sum_{w\in W^{(2)}(M_{II})}{ \hskip-10pt \det(w)
\left(\exp(-\pi i (w(\rho),z))
- \hskip-10pt  \sum_{ a \in M_{II}\cap \br_{++}\M_{II}}
{m(a)\exp(-\pi i (w(\rho+a),z))}\right)},
$$
where $z=z_3f_2+z_2f_3+z_1f_{-2} \in M_0 \otimes \br +i V^+(M_0)=
M_{II} \otimes \br +i V^+(M_{II})$ and
$m(a) \in \bz$.

{\rm {(b)}} For any primitive $a_0 \in M_{II}\cap \br_{++}\M_{II}$
with $(a_0,a_0)=0$, we have the identity
$$
1- \sum_{t \in \bn}{m(ta_0)q^t}=
\prod_{k \in \bn}{(1-q^k)^{\tau(ka_0)}}
$$
of power series of one variable $q$, where all $\tau(ka_0)\in \bz$.
More exactly,
$$
\tau(a)=9\ \ \text{for any}\ \  a \in M_{II}\cap \br_{++}\M_{II}\
\text{with}
\ (a,a)=0.
$$
\endproclaim

\demo{Proof} Statement (a). The lattice
$M_{II}^\ast=\bz (f_2/2)+\bz (f_3/2) +\bz (f_{-2}/2)=\frac{1}{2}M_0$.
Thus, for $n,l,m \in \bz$, we can rewrite
$$
\split
&(1/64)f(n,l,m)\exp(\pi i (nz_1+lz_2+mz_3))=\\
=&(1/64)f(n,l,m)\exp(-\pi i (nf_2-lf_3/2+mf_{-2},z))
=-m(a)\exp(-\pi i (\rho + a, z)),
\endsplit
$$
where $a=(n-1)f_2-(l-1)f_3/2+(m-1)f_2 \in M_0^\ast =\frac{1}{2}M_{II}$
and $m(a)=-(1/64)f(n,l,m)$. By \thetag{1.6} and
\thetag{1.7}, we have that
$\rho +a \in V^+(M_0)$, $m(a) \in \bz$, and
$m(0)=-1$.  By Proposition 2.1, (c), we get
$$
{\tsize\frac1{64}}\Delta_5(Z)=\sum_{w\in W^{(2)}(M_{II})}{\det(w)
\left(
-\sum_{ \rho +a \in M_{II}^\ast \cap \br_+\M_{II}}
{m(a)\exp(-\pi i (w(\rho+a),z))}\right)}.
$$
Since for this sum  $\rho +a \in M_{II}^\ast \cap \br_+\M_{II}$,
we have that $(\rho+a, \delta_i)\le 0$, $i=1,2,3$. If
$(\rho+a, \delta_i)=0$, then the corresponding Fourier coefficient
$m(a)=0$, since $\Delta_5(Z)$ is anti-invariant with respect to
$s_{\delta_i}$. Thus, we can suppose that
$(\rho+a, \delta_i)$ is integral and $(\rho+a, \delta_i)<0$,
$i=1,2,3$. By definition of $\rho$, we
then get that $(a, \delta_i)\le 0$. It follows that $a \in \br_+Q_+^\ast$.
By \thetag{2.2}, then $a \in \br_+Q_+^\ast\cap \overline{V^+(M_0)}=
\br_+\M_{II}$. It follows that $a \in \br_{++}\M_{II}$ if $a \not=0$.
If $a=0$, we have $m(a)=-1$.
By the congruences $m,n,l \equiv 1\mod 2$ of \thetag{1.6}, we have
$a \in 2M_0^\ast=M_{II}$.

Statement (b). Primitive elements
$a_0 \in M_{II}\cap \br_{++}\M_{II}$ with
$(a_0,a_0)=0$ correspond to three infinite
vertices of $\M_{II}$. By Statement (c) of Proposition 2.2, the
group $A(P(\M_{II})_{\pr})$ is transitive on these three vertices and
the corresponding primitive elements $a_0$ and
preserves the Fourier expansion (a). Thus, it is sufficient to prove
the identity (b) for $a_0=2f_2$ which is one of these three primitive
elements of $M_{II} \cap \br_{++}\M_{II}$ (all these three primitive
elements are $2f_2$, $2f_{-2}$, $2f_2-2f_3+2f_{-2}$).
For $a_0=2f_2$, Statement (b) is equivalent to \thetag{1.8}.
It follows Statement (b).
\enddemo

\head
\S~3. The Kac--Moody algebra $\geg(M_{II}, P(\M_{II}))$ and
its $\Delta_5(Z)$-correction $\geg=\geg (M_{II}, {}_s\Delta)$
where ${}_s\Delta=P(\M_{II})\cup {}_s\Delta^{\im}$.
\endhead

The matrix \thetag{2.1} is the Gram matrix of elements $P(\M_{II})$, and
we denote this Gram matrix as
$G(P(\M_{II}))$. This is a symmetric generalized Cartan matrix
(i.e. it is integral, has only $2$ on the diagonal and
only non-positive integers out of the diagonal).
This matrix defines
the corresponding Kac--Moody algebra
$\geg (M_{II},P(\M_{II}))=:\geg (G(P(\M_{II}))$
(see \cite{Ka1} and  \cite{Mo}).
We want to ``correct" this algebra using Theorem 2.3.
Some examples of this type constructions
were first found by R. Borcherds
(see \cite{Bo3}, \cite{Bo5}, \cite{Bo6}, \cite{Bo7}).
Its formulation in an appropriate general setting see in \cite{N5}.
This construction uses
generalized Kac--Moody algebras which were introduced and studied
by R. Borcherds in \cite{Bo1} (see also \cite{Bo2}).
For our situation we need the corresponding superalgebra analog of this
construction. It is a combination of Borcherd's
construction and of V. Kac (\cite{Ka2}, \cite{Ka3}).
We give this construction for our case below and refer
to \S~6, Appendix for proofs.

\vskip5pt

Using coefficients $m(a)$ and $\tau(a)$ of Theorem 2.3,
we introduce (where the index $s$ means ``simple") the sets
$$
\split
{}_s \Delta_{\overline 0}^{\im}=&\{m(a)a\ |\
a\in M_{II}\cap \br_{++}\M_{II},\ (a,a)<0\ \text{and}\  m(a)>0\}\cup\\
\cup &\{\tau(a)a\ |\
a\in M_{II}\cap \br_{++}\M_{II},\ (a,a)=0\ \text{and}\  \tau (a)>0\}
\endsplit
$$
where $ka$ for $k\in \bn$ means that we repeat $a$ exactly $k$ times, and
$$
\split
{}_s\Delta_{\overline 1}^{\im}=&
\{m(a)a\ |\
a\in M_{II}\cap \br_{++}\M_{II},\ (a,a)<0\ \text{and}\  m(a)<0\}\cup \\
\cup &\{\tau(a)a\ |\
a\in M_{II}\cap \br_{++}\M_{II},\ (a,a)=0\ \text{and}\  \tau (a)<0\}
\endsplit
$$
where $ka$ for $-k\in \bn$ means that we repeat the element $a$ exactly
$-k$ times. The minus sign of $k$
here denote that all this $-k$ elements $a$
are considered to be odd.
The elements of ${}_s\Delta_{\overline 0}^{\im}$ and
${}_s\Delta_{\overline 1}^{\im}$ are considered to be even
and odd respectively. They are
called {\it even} and {\it odd imaginary simple roots} respectively.
We denote
$$
{}_s\Delta^{\im}={}_s\Delta_{\overline 0}^{\im}\cup
{}_s\Delta_{\overline 1}^{\im}.
$$
Elements of ${}_s\Delta^{\im}$ are called
{\it imaginary simple roots}.
To unify notations we denote
$_s\Delta_{\overline 0}^{\re}={}_s\Delta^{\re}=P(\M_{II})$.
Elements of ${}_s\Delta^{\re}$ are called
{\it real simple roots}.
They are considered to be even (for our situation). We denote as
${}_s\Delta ={}_s\Delta^{\re} \cup {}_s\Delta^{\im}$
the whole set of simple roots. Thus, simple roots are either real
even or imaginary even or imaginary odd.
By the construction, elements of ${}_s\Delta$ give the
corresponding elements of
$M_{II}\subset M_{II}\otimes \br$. For the further
construction of the corresponding Kac--Moody superalgebra
we need that $(r,r)>0$ if $r \in {}_s\Delta^{\re}$,
$(r,r) \le 0$ if $r \in {}_s\Delta^{\im}$, and
$(r,r^\prime)\le 0$ for all different
$r,r^\prime \in {}_s\Delta$.
This is true by the construction.
Moreover, for
$r, r^\prime \in {}_s\Delta$ we need
$$
2(r,r^\prime)/(r,r)\in \bz \quad \ \text{if\ }(r,r)>0.
$$
This is valid because here $(r,r)=2$ and
$r^\prime \in M_{II}$ where $M_{II}$ is generated by
elements of ${}_s\Delta^{\re}$
(which are $\delta_1,\delta_2,\delta_3$).
The generalized Kac--Moody superalgebra without odd real simple roots
$\geg=\geg^{\prime\prime}(M_{II}, {}_s\Delta)$
(which is constructing) is a Lie superalgebra generated by
$h_r, e_r, f_r$ where $r \in {}_s\Delta$. Here all $h_r$ are even,
$e_r, f_r$ are even (respectively odd) if
$r$ is even (respectively odd).
These elements have the following defining relations:

\vskip5pt

(1) The map $r \mapsto h_r$ for $r\in {}_s\Delta$ gives an embedding
of $M_{II}\otimes \br$ to $\geg^{\prime\prime}(M_{II}, {}_s\Delta)$ as
an abelian subalgebra (it is even since all $h_r$ are even).
In particular, all elements $h_r$ commute.

(2) $[h_r, e_{r^\prime}]=(r, r^\prime)e_{r^\prime}$, and
$[h_r, f_{r^\prime}]=-(r,{r^\prime})f_{r^\prime}$.

(3) $[e_r, f_{r^\prime}]=h_r$ if $r=r^\prime$, and is $0$ if
$r \not=r^\prime$.

(4) $(\text{ad\ } e_r)^{1-2(r,r^\prime)/(r,r)}e_{r^\prime }=
(\text{ad\ } f_r)^{1-2(r,r^\prime)/(r,r)}f_{r^\prime }=0\
\text{if $r\not= r^\prime$ and $(r,r)>0$}$ \hfil\hfil
\newline
\phantom{(4)(4)} (equivalently, $r \in {}_s\Delta^{\re}$).

(5) If $(r,r^\prime)=0$, then $[e_r, e_{r^\prime}]=[f_r,f_{r^\prime}]=0$.

\vskip5pt

The superalgebra $\geg=\geg ^{\prime\prime}(M_{II}, {}_s\Delta)$
is graded by $M_{II}$ as follows. Let
$$
\widetilde{Q}_+ =\sum_{\alpha \in {}_s\Delta}{\bz_+\alpha}\subset M_{II}
$$
be the integral cone (semi-group) generated by all simple roots.
(For our case $\widetilde{Q}_+$ coincides with the integral
cone  $Q_+=\bz_+{}_s\Delta^{\re}$ generated by real simple roots
${}_s\Delta^{\re}=\{\delta_1, \delta_2, \delta_3\}$,
compare with \thetag{2.2}.) We have
$$
\geg=\left( \bigoplus_{\alpha \in \widetilde{Q}_+} {\geg_\alpha} \right)
\bigoplus \left({M_{II}}\otimes \br \right) \bigoplus
\left( \bigoplus_{\alpha \in -\widetilde{Q}_+} {\geg_\alpha} \right)
$$
where $e_r$ and $f_r$ have degree $r\in \widetilde{Q}_+$ and
$-r \in -\widetilde{Q}_+$ respectively, $r \in {}_s\Delta$;
and $\geg_0=M_{II}\otimes \br$.
The $0\not=\alpha \in \pm \widetilde{Q}_+$ is called  a
{\it root} if $\geg_\alpha$ is
non-zero. Let $\Delta$ be the set of all roots and
$\Delta_{\pm}=\Delta\cap \pm \widetilde Q_+$.
For a root $\alpha \in \Delta$ we denote
$\mult_\0o \alpha =\dim \geg_{\alpha,\0o}$,
$\mult_\1o \alpha = -\dim \geg_{\alpha,\1o}$ and
$$
\mult~\alpha=\mult_\0o \alpha +\mult_\1o \alpha=
\dim \geg_{\alpha ,\0o} - \dim \geg_{\alpha,\1o}.
\tag3.1
$$
By our construction and
Statement 6.8$^\prime$ in \S~6, it follows the
denominator identity for $\geg$: $\Phi (z)=:$
$$
\multline
\hskip-10pt \sum_{w\in W^{(2)}(M_{II})}{\hskip-15pt \det(w)
\left(\exp(-2\pi i (w(\rho),z))
- \hskip-10pt \sum_{ a \in M_{II} \cap \br_{++}\M_{II}}\hskip-10pt
{m(a)\exp(-2\pi i (w(\rho+a),z))}\right)}\\
=\exp{\left(-2\pi i(\rho,z)\right)}
\prod_{\alpha \in \Delta_+}
{\left( 1-\exp{ \left(-2\pi i (\alpha , z)\right)}\right)^{\mult~\alpha}},
\endmultline
\tag3.2
$$
where $z \in \Omega(M_{II})=M_{II}\otimes \br +iV^+(M_{II})$.
The function $\Phi (z)$ is called the {\it denominator function}.
Thus, using these considerations and Theorem 2.3, we get

\proclaim{Proposition 3.1} We have the following equality between
the $\Delta_5(Z)$ and the denominator function $\Phi (z)$
of the correcting generalized Kac--Moody superalgebra \linebreak
$\geg=\geg ^{\prime\prime} (M_{II}, {}_s\Delta)$:
$$
{\tsize\frac{1}{64}} \Delta_5(2Z)=\Phi (z).
$$
Thus, the denominator function $\Phi (z)$ of the correcting
generalized Kac--Moody superalgebra
$\geg ^{\prime\prime} (M_{II}, {}_s\Delta)$ is a Siegel modular form.
On the other hand, the classical Siegel modular form
$\Delta_5(Z)$ is the denominator function of a
generalized Kac--Moody superalgebra.
\endproclaim

\remark{Remark 3.2} The denominator function $\Phi(z)$ is defined
in the complexified cone \linebreak $\Omega (V^+(M_{II}))$ (we
remark that to have this property we again need \thetag{2.2}
and also existence of a lattice Weyl vector
$\rho \in M_{II}\otimes \bq$).
This cone does not have
a canonical embedding as a cusp to a IV type domain. This embedding is
defined up to changing $z$ to $tz$ where $t \in \bn$. This is an
explanation of the coefficient $2$ for the equality of
Proposition 3.1.  See \cite{N5} for the appropriate general setting.
\endremark

\vskip5pt

The equality \thetag{3.2} shows that there exists a product formula for
$\Delta_5(Z)$
$$\multline
\frac{1}{64} \Delta_5(Z) \hskip-1pt =
\\
\hskip-10pt
\sum_{w\in W^{(2)}(M_{II})}{ \hskip-10pt \det(w)
\left(\exp(-\pi i (w(\rho),z))
- \hskip-10pt \sum_{ \hskip-10pt a \in M_{II} \cap \br_{++}\M_{II}}
{ \hskip-10pt m(a)\exp(-\pi i (w(\rho+a),z))}\right)}\\
=\exp(-\pi i (\rho, z)) \prod_{\alpha \in \Delta_+}
{\left( 1-\exp({-\pi i (\alpha , z)})\right)^
{\mult~\alpha}}.
\endmultline
\tag3.3
$$
On the other hand, using automorphic properties of
$\Delta_5(Z)$, we may have a hope to calculate this product formula and
calculate multiplicities $\mult~\alpha$ for products
\thetag{3.2} and \thetag{3.3}. We do it in \S~4.  We will show
that for $\alpha=2nf_2-lf_3+2mf_{-2} \in M_{II}$ (i.e. $n,m,l \in \bz$),
$$
\mult~\alpha=f(nm,l)
\tag3.4
$$
where $f(k_1,k_2)$ are Fourier coefficients of the weak Jacobi function
$\phi_{0,1}(z_1,z_2)$ of weight $0$ and index $1$. See \thetag{4.11}
and Theorem 4.1.

\head
\S~4.  The product formula for $\Delta_5(Z)$ and root spaces
multiplicities $\mult~\alpha$ for by
$\Delta_5(Z)$ corrected Kac--Moody algebra $\goth g$
\endhead

The modular form $\Delta_5(Z)$ has the following Fourier-Jacobi expansion
$$
\Delta_5\bigr(\pmatrix z_1&z_2\\z_2&z_3\endpmatrix\bigl)=
\sum\Sb m>0\\ m\equiv 1\,mod\, 2\endSb
\phi_{5,m}(z_1,\,z_2)\exp{(\pi i mz_3)}.
\tag4.1
$$
This is the   Fourier expansion  along the one dimensional cusp
of domain $Sp_4(\bz)\setminus \bh$ corresponding to the following  maximal
parabolic subgroup
$$
\gi=
\left\{\pmatrix *&0&*&*\\
                *&*&*&*\\
                *&0&*&*\\
                0&0&0&* \endpmatrix
\in Sp_4(\Bbb Z)\right\}.
$$
The group $\Gamma^J=\gi/\{\pm E_4\} \cong SL_2(\bz)\ltimes H(\bz)$
is called  Jacobi group.
We realize $SL_2(\bz)$ as a subgroup of $Sp_4(\bz)$
throught the embedding
$$
*:\,\pmatrix a&b\\c&d\endpmatrix\mapsto
\pmatrix a&b\\c&d\endpmatrix^*=
\pmatrix a&0&b&0\\
               0&1&0&0\\
                c&0&d&0\\
                0&0&0&1 \endpmatrix.
\tag4.2
$$
$H(\bz)$ is the integral  Heisenberg group
$$
H(\Bbb Z)\cong\left\{h(p,q;r)=\pmatrix 1&0&0&q\\
               p&1&q&r\\
                0&0&1&-p\\
                0&0&0&1 \endpmatrix \in Sp_4(\Bbb Z)\right\}.
\tag4.3
$$
The restriction of multiplier system $v$ on $\gi$ define a character
$$
v_\infty:SL_2(\bz)\ltimes H(\bz)\to \{\pm 1\}.
$$
The functions $\phi_{5,m}(z_1,\,z_2)$ are Jacobi forms of
{\bf half-integral index}.
They have the following invariant properties, which one can take
as a definition of Jacobi forms of index $\dsize\frac{m}2$:
$$
\align
\kern-30pt
\phi_{5,m}(z_1,\,z_2)\exp{(\pi i\, mz_3)}|_5\, h(p,q;r)&=
(-1)^{p+q+pq+r}\phi_{5,m}(z_1,\,z_2)\exp{(\pi i \,mz_3)}
\tag4.4\\
\kern-30pt
\phi_{5,m}(z_1,\,z_2)\exp{(\pi i\, mz_3)}|_5\, g^*&=
v( g^*)\phi_{5,m}(z_1,\,z_2)\exp{(\pi i\, mz_3)}
\tag4.5
\endalign
$$
for $ h(p,q;r)\in H(\bz)$ and $ g\in SL_2(\bz)$.
We set
$$
F|_k\,M\,(Z):=\hbox{det\,}(CZ+D)^{-k}\,F\bigl((AZ+B)(CZ+D)^{-1}\bigr)
$$
for any
$M=\left(\smallmatrix A&B\\C&D\endsmallmatrix\right)\in Sp_4(\br)$
and any function
$F:\,\Bbb H_2\to \bc$.
One can rewrite the last two equations  in the following form
$$
\align
\phi_{5,m}(\frac{az_1+b}{cz_1+d},\,\frac{z_2}{cz_1+d})&=
v(g^*)(cz_1+d)^{5}
\exp{(\pi i\,\frac{ cmz_2^2}{cz_1+d})}\,\phi_{5,m}(z_1 ,\,z_2 ),
\\
\vspace{2\jot}
\phi_{5,m}(z_1, z_2+pz_1+q)&=
(-1)^{p+q+pq(m+1)}\exp{(-\pi i m(p^2z_1 +2p z_2))}\,\phi_{5,m}(z_1, z_2),
\endalign
$$
where $p,q\in \bz$ and
$g=\left(\smallmatrix a&b\\c&d\endsmallmatrix\right) \in SL_2(\bz)$.
One can easy check that
$$
v(\pmatrix 1&b\\0&1\endpmatrix^*)=(-1)^b
\qquad\text{and}\qquad
v(\pmatrix 0&-1\\1&0\endpmatrix^*)=-1,
$$
since
$$
\pmatrix 0&-1\\1&0 \endpmatrix^*=
\biggl[\pmatrix 0&E_2\\-E_2&0\endpmatrix
\pmatrix E_2&B_1\\0&E_2\endpmatrix\biggr]^3
\pmatrix 0&-E_2\\E_2&0\endpmatrix
\qquad B_1=\pmatrix 1&0\\0&0\endpmatrix.
$$

The first Fourier-Jacobi coefficient $\phi_{5,1}(z_1,\,z_2)$
is connected with the denominator formula for a
generalized affine Kac--Moody algebra.
In order to get this connection let us take
the following Jacobi theta-series
$$
\vartheta_{11}(z_1,z_2)=\sum\Sb n\in \bz\endSb
\,(-1)^{n}
\exp{\bigl(\frac{\pi i\, }{4}(2n+1)^2 z_1+\pi i\,(2n+1)z_2\bigr)}.
$$

We recall  a well known variant of the  Jacobi triple formula:
$$
\prod_{n\ge 1}(1-q^{n-1}r)(1-q^n r^{-1})(1-q^n)=
\sum_{m\in \bz}(-1)^m q^{\frac 1{2}m(m-1)}r^m.
$$
 From that the product formula for $\vartheta_{11}$ follows
$$
\vartheta_{11}(z_1,\,z_2)=
-q^{1/8}r^{-1/2}
\prod_{n\ge 1}\,(1-q^{n-1} r)(1-q^n r^{-1})(1-q^n).
\tag4.6
$$
It is easy to see from definition that $\vartheta_{11}(z_1,z_2)$
and $\phi_{5,1}(z_1,z_2)$
have  the same functional equation with respect to $h(p,q;r)$:
$$
\vartheta_{11}
(z_1,z_2+pz_1+q)=
(-1)^{p+q}\exp{(-\pi i\,(p^2z_1 +2p z_2))}\,\vartheta_{11}(z_1, z_2).
$$
Moreover $\vartheta_{11}$ satisfies
the following transformations formulae with respect to
\newline
generators
$\left(\smallmatrix 0&-1\\1&0\endsmallmatrix\right)$ and
$\left(\smallmatrix 1&1\\0&1\endsmallmatrix\right)$ of $SL_2(\bz)$:
$$
\align
\vartheta_{11}(-\frac 1{z_1},\, \frac {z_2}{z_1})&=
-i\sqrt{-iz_1}\,\exp{(\pi i\,\frac {z_2^2}{z_1})}\,
\vartheta_{11}(z_1,z_2),
\\
\vartheta_{11}(z_1+1,\,z_2)&
=\exp{(\frac {\pi i}{4})}\,\vartheta_{11}(z_1,\,z_2).
\endalign
$$
The Dedekind eta-function
$$
\eta(\tau)=
\exp{(\frac{\pi i \tau}{12})}\,\prod_{n\ge 1}\ (1-\exp{(2\pi i\,n\tau)})
$$
satisfies the similar functional equations
$$
\eta(-\frac 1{\tau})=-i\sqrt{-i\tau}\,\eta(\tau),\qquad
\eta(\tau+1)=\exp{(\frac {\pi i}{12})}\,\eta(\tau).
$$
Thus $\phi_{5,1}(z_1,\,z_2)$ and
$$
\psi_{5,\frac 1{2}}(z_1,z_2)
=\eta(z_1)^9\,\vartheta_{11}(z_1,\,z_2)
\tag4.7
$$
are   Jacobi cusp forms of index $\dsize\frac {1}2$ with
the same character $v_\infty:\,\Gamma^J\to\{\pm 1\}$.
The squares of these Jacobi forms
are Jacobi cusp forms of weight $10$ and index $1$.
Up to a constant there is only one such  a  form.
This is the  first Fourier-Jacobi coefficient
$\phi_{10,1}(z_1,\,z_2)$ of
$F_{10}(Z)$ from \thetag{1.1}  (see \cite{EZ}).
Comparing a pair of their Fourier coefficients
we  obtain
$$
{\tsize\frac{1}{64}}\phi_{5,1}(z_1,\,z_2)=\psi_{5,\frac 1{2}}(z_1,z_2)=
-q^{1/2}r^{-1/2}
\prod_{n\ge 1}\,(1-q^{n-1} r)(1-q^n r^{-1})(1-q^n)^{10},
\tag4.8
$$
where
$$
q=\exp{(2\pi i\, z_1)}, \qquad r=\exp{(2\pi i \,z_2)}.
$$
If we write down the right hand side as a formal series in $r$
then the $q$-sum in \thetag{1.8} is the coefficient of $r^{1/2}$.
Using the Jacobi triple formula we get \thetag{1.8}.

Let $\phi_{12,1}(z_1, z_2)$ be a  Jacobi cusp form of weight $12$
and index one (the first Fourier-Jacobi coefficient of $F_{12}(Z)$ from
\thetag{1.1}).
It is known that up to a constant this is the only one such a form.
Its Fourier coefficient can be calculated according to
the following formula (see \cite{EZ})
$$
\phi_{12,1}(z_1, z_2)=\frac 1{144}
\bigr(E_4^2(z_1)E_{4,1}(z_1,z_2)-E_6(z_1)E_{6,1}(z_1,z_2)\bigl),
\tag4.9
$$
where
$$
E_4(z_1)=1+240 \sum_{n\ge 1}\sigma_3(n)q^n, \qquad
E_6(z_1)=1-504 \sum_{n\ge 1}\sigma_5(n)q^n
\tag4.10
$$
are  Eisenstein series for $SL_2(\bz)$
and $E_{k,1}(z_1,z_2)$ is  Eisenstein series with integral
coefficients
of  weight $k$  and index one
for the Jacobi group
$$
E_{k,1}(z_1,z_2)=\sum\Sb n,l\in \bz\\ 4n-l^2\ge 0\endSb
\frac {H(k-1,4n-r^2)}{\zeta(3-2k)} \exp{(2\pi i \,(nz_1+lz_2))},
$$
where $H(k,N)=L_{N}(1-k)$ are H. Cohen's numbers (see \cite{EZ}).
The form $\phi_{12,1}$ has integral and coprime coefficients:
$$
\phi_{12,1}(z_1, z_2)=(r^{-1}+10+r)q +
(10r^{-2}-88r^{-1}-132-88r+10r^{-2})q^2 +\dots\ .
$$
Let us introduce a function with integral Fourier coefficients
$$
\phi_{0,1}(z_1, z_2):=\phi_{12,1}(z_1, z_2)/\Delta_{12}(z_1)=
\sum_{n\ge 0,\, l\in \Bbb Z} f(n,l)\,\exp{(2\pi i \,(nz_1+lz_2))},
\tag4.11
$$
where
$$
\Delta_{12}(z_1)=q\prod_{n\ge 1}(1-q^n)^{24}
$$
is the $SL_2(\Bbb Z)$-cusp form of weight $12$.
$\phi_{0,1}(z_1, z_2)$ is a weak Jacobi function of weight $0$
and index $1$.
This function satisfies the same functional
equation as holomorphic Jacobi forms and has  nonzero Fourier coefficients
only with indexes $(n,l)\in \bz$ such that
 $n\ge 0$ (since $\phi_{12,1}$ is a cusp form)
and $4n-l^2\ge -1$.
The weight is even, thus   $f(n,l)=f(n,-l)$
and $f(n,l)$ depends only on $4n-l^2$.
Moreover
$$
\phi_{0,1}(z_1,z_2)=(r^{-1}+10+r)+q(10r^{-2}-
64r^{-1}+108-64r+10r^2)+\dots\ .
$$
All facts from the theory of Jacobi forms mentioned above
 one can find in \cite{EZ}.

\proclaim{Theorem 4.1}The following formula is valid:
$$
{\tsize\frac 1{64}}\,\Delta_5(Z)=\exp{(\pi i\,(z_1+z_2+z_3))}
\kern -10pt
\prod
\Sb n,\,l,\, m\in \Bbb Z\\
\vspace{0.5\jot}
(n,l,m)>0\endSb
\bigl(1-\exp{(2\pi i\,(nz_1+lz_2+mz_3))}\bigr)^{f(nm,l)},
$$
where $(n,l,m)>0$ means that $n\ge 0$, $m\ge 0$,
$l$ is an  arbitrary integral if $n>0$ or $m>0$ and
$l<0$ if $n=m=0$.
\endproclaim

\remark{Remark 4.2}The condition $(n,l,m)>0$ means that the product
is taken over the set of  positive roots $\Delta_+$
(see \thetag{3.2}).
\endremark

\demo{Proof of Theorem} The holomorphic Jacobi form  has the following
Fourier expansion
$$
\align
\phi_{12,1}(z_1,z_2)&=\sum_{\mu\,mod\,2}
\sum\Sb 4n-l^2>0\\ l\equiv \mu\,mod 2\endSb
c_{\mu}(4n-l^2)\exp{(2\pi i\,(n z_1+l z_2))}\\
\vspace{2\jot}
{}&=
h_0(z_1)\vartheta_{1,0}(z_1, z_2)+
h_1(z_1)\vartheta_{1,1}(z_1, z_2),
\endalign
$$
where $h_\mu(z_1)=\sum_{m>0} c_{\mu}(m)\exp{(\frac{\pi i}2\, m z_1)}$
are  modular forms of weight $\frac{23}2$
and $\vartheta_{1,\mu}(z_1, z_2)$ are the standard Jacobi theta-functions.
Thus the Fourier coefficients of
$\phi_{0,1}(z_1,z_2)$
have the  following asymptotic
$$
f(n,l)= O(e^{\sqrt{4n-l^2}}).
$$
Using this estimate, \thetag{2.2} and arguments of
\cite{Ka1, \S~10.6},
we prove that the product of Theorem 4.1 absolutely converges on any
neighborhood of the zero-dimensional  cusp $\infty$ of $Sp_4(\bz)$.

Let us denote the  product in the formula of Theorem 4.1
 by $P(z_1, z_2, z_3)$ and decompose
it in two factors corresponding to the cases $m=0$ and $m>0$
$$
\multline
P(z_1,z_2,z_3)=
\exp{(\pi i\,(z_1+z_2+z_3))}\,
\prod\Sb n> 0,\, l\in \Bbb Z\\
\vspace{0.5\jot}n=0,\,l<0\endSb
\big(1-\exp{(2\pi i\,(nz_1+lz_2))}\bigr)^{f(0,l)}\\
\vspace{2\jot}
\times
\prod_{n\ge 0,\,m>0,\, l\in \Bbb Z}\,
\bigl(1-\exp{(2\pi i\,(nz_1+lz_2+mz_3))}\bigr)^{f(nm,l)},
\endmultline
\tag4.12
$$
where $f(0,0)=10$, $f(0,-1)=1$ and $f(0,l)=0$
if $l<-1$.

Let us introduce the following Hecke operators
$$
T_{-}(m)=\sum\Sb\alpha\beta=m\\ \alpha\,|\,\beta\endSb
\Gamma_\infty
\pmatrix
\alpha&0&0&0\\0&m&0&0\\0&0&\beta&0\\0&0&0&1
\endpmatrix\Gamma_\infty=
\sum\Sb ad=m\\
\vspace{0.5\jot} b\,mod\,d\endSb \Gamma_\infty
\pmatrix
a&0&b&0\\0&m&0&0\\0&0&d&0\\0&0&0&1
\endpmatrix,
$$
which are the ``minus''-embedding of the usual Hecke operators
$T(m)$ for $GL_2(\bz)$ (see \cite{G3}).
By the definition
$$
\pmatrix
a&0&b&0\\0&m&0&0\\0&0&d&0\\0&0&0&1
\endpmatrix
<\pmatrix z_1&z_2\\z_2&z_3\endpmatrix>=
\pmatrix
\frac{a z_1+b}d&a z_2\\ a z_2&mz_3
\endpmatrix,
$$
therefore for each Jacobi form
$\widetilde\phi(Z)=\phi(z_1,z_2)\exp{(2\pi i\,tz_3)}$ of
weight $k$ and index $t$ ($t\in \bq$) the   function
$$
\bigl(\widetilde \phi\,|_k\, \,T_{-}(m)\bigr)(Z)=
m^{2k-3} \sum\Sb ad=m\\ \vspace{0.5\jot} b\,mod\,d\endSb
d^{-k}\,
\phi (\frac{az_1+b}d,\, az_2)
\exp{(2\pi i\, mt z_3)}
$$
is a Jacobi form of index $mt$. (We put a normalizing factor
 $m^{2k-3}$
in the definition of $|_k\, T_-(m)$ like in \cite{G3}).

Let us calculate the Fourier expansion of the logarithm
of the second factor from \thetag{4.12}:
$$
\multline
\hbox{log\,}\big(
\prod\Sb n\ge 0,\,m> 0,\,  l\in \Bbb Z\endSb
\dots\bigr)=
-\kern -8pt
\sum\Sb n\ge 0,\,m> 0\\
\vspace{1\jot}
l\in \Bbb Z\endSb
f(nm,l)\ \sum_{e\ge 1}\frac 1{e}\exp{(2\pi i e\,(nz_1+lz_2+mz_3))}
\\
= - \kern -5pt
\sum\Sb a\ge 0,\,c> 0\\
\vspace{1\jot}
b\in \Bbb Z\endSb
\ \sum\Sb t\,|\,(a,b,c)\endSb
t^{-1}f(\frac{ac}{t^2},\,\frac{b}{t})
\exp{(2\pi i\,(az_1+bz_2+cz_3))}.
\endmultline
$$
The last sum can be written as a sum of Jacobi forms using
the operators $T_{-}(m)$.
By definition
$$
\align
m^2\,\bigl(\widetilde{\phi}_{0,1}\,|_0\,T_-(m)\bigr)(Z)
&= \sum_{ad=m}m^{-1}
\hskip-10pt\sum\Sb n,\,l\\ \vspace{0.5\jot}
b\mod d\endSb
f(n,l)\,\exp{(2\pi i\,(n\,\frac{az_1+b}{d}+la z_2+mz_3))}\\
\vspace{2\jot}
{}&=
\sum_{ad=m}a^{-1}\sum_{n,l}\,f(dn,l)
\exp{(2\pi i (an z_1+al z_2+ ad z_3))}.
\endalign
$$
 From that one easily gets
$$
\hbox{log\,}\big(
\prod\Sb n\ge 0,\,m> 0,\, l\in \Bbb Z\endSb
\dots\bigr)=-\sum_{m\ge 1}m^2\,
\bigl(\widetilde{\phi}_{0,1}\,|_0\,T_-(m)\bigr)(Z).
$$
This expansion shows us that  the second factor in \thetag{4.12}
is invariant with respect to the action of the elements of
$\Gamma_\infty$.

The first factor in \thetag{4.12} is equal to
$\widetilde{\psi}_{5,\frac1{2}}(Z)$ (see \thetag{4.8}).
Therefore
$P(z_1,z_2,z_3)$ transforms like a $\Gamma^J$-modular form of
weight $5$ with a character  $v_\infty:\Gamma^J\to \{\pm 1\}$.

There is  the only factor in $P(z_1,z_2,z_3)$ with $n=m=0$ and
we may rewrite the product in the following form
$$
\multline
P(z_1,z_2, z_3)=
\big(\exp{(\pi i\,(z_1+z_2+z_3))}
-\exp{(\pi i\,(z_1-z_2+z_3))}\bigr)\times
\\
\prod \Sb n,m\ge 0,\, l\in \Bbb Z\\
\vspace{1\jot}
n>0\vee m>0\endSb
\bigl(1-\exp{(2\pi i\,(nz_1+lz_2+mz_3))}\bigr)^{f(nm,l)}.
\endmultline
$$
Thus
$
P(z_1,z_2, z_3)=P(z_3,z_2, z_1)
$.
The element
$$
I=\pmatrix U&0\\ 0&U\endpmatrix,\quad\text{where }\
U=\pmatrix 0&1\\ 1&0\endpmatrix,
$$
realizes  this  change of the variables.
Hence
$$
P(z_1,z_2, z_3)\,|_5\, I=-P(z_1,z_2, z_3)\qquad\text{and}\qquad v_P(I)=1.
$$
$\Gamma^J$ and $I$ generate the group $PSp_4(\bz)$. Therefore
$P(z_1,z_2, z_3)$ is a Siegel modular form of weight $5$ with the same
multiplier system as $\Delta_5(Z)$. Comparing the first
Fourier coefficients we finish the proof of the theorem.
\enddemo

\medskip
There is a better way to define the  Fourier expansion of $\Delta_5(Z)$,
than the definition  \thetag{1.2}.  This function belongs to
the so-called Maass subspace (see \cite{M3}). It means that
$\frac 1{64}\Delta_5(Z)$ can be written as a lifting
of its ``first'' (it is better to say ``one half'')
Fourier-Jacobi coefficient
$$
\psi_{5,\frac 1{2}}(z_1,z_2)=\eta(z_1)^9\,\vartheta_{11}(z_1,z_2)=
\sum_{n,l\,\equiv 1\,mod\  2} g(n,l)\,\exp{(\pi i\,(nz_1+lz_2))}.
$$
More exactly
$$
\gather
{\tsize\frac 1{64}}\Delta_5(Z)=\sum\Sb n,l,m\equiv 1\, mod\  2\\
                    n,m>0,\ 4mn-l^2>0\endSb
\sum_{d|(n,l,m)} d^4\, g(\frac{nm}{d^2},\,\frac{l}{d})
\exp{(\pi i\,( nz_1+lz_2+mz_3))}=
\\
\vspace{2\jot}
\kern-10pt\exp{(\pi i\,(z_1+z_2+z_3))}\,
\prod
\Sb n,\,l,\,m\in \Bbb Z\\
\vspace{0.5\jot}
(n,l,m)>0\,\endSb
\bigl(1-\exp{(2\pi i\,(nz_1+lz_2+mz_3))}\bigr)^{f(nm,l)}
\tag4.13
\endgather
$$
where $f(\cdot,\cdot)$ are the Fourier coefficients
of the weak Jacobi form $\phi_{0,1}$.

\head
\S~5. The second example of Siegel automorphic form correction of a
Lorentzian Kac--Moody algebra
\endhead

Here we want to describe briefly another example of Siegel automorphic
form correction of a Lorentzian Kac--Moody algebra. Actually,
the case we have considered before and this example belong to one
series which we hope to describe and study carefully in the
corresponding publication later. The main difference of this example is
that Siegel modular form we will use for this correction is not known.
To construct this modular form we use the arithmetic lifting of
Jacobi forms which was developed in \cite{G1} -- \cite{G3}.
This lifting is the advanced generalization of Maass's lifting.
This lifting is valid for an arbitrary domain of type 4,
and we hope to apply this lifting for other cases later.

\vskip10pt

We start with the description of the Lorentzian Kac--Moody
algebra we want to correct, its real roots and Weyl group (i.e. the
corresponding reflection group in hyperbolic plane).

We consider
an even unimodular hyperbolic plane $\bz f_2 \oplus \bz f_{-2} \cong U$
with
$(f_2,f_2)=(f_{-2},f_{-2})=0$, $(f_2,f_{-2})=-1$,
and a one-dimensional
lattice
$\bz f_3$ with $(f_3, f_3)=4$, and their orthogonal sum
$M_0=\bz f_2 \oplus \bz f_3 \oplus \bz f_{-2}$.
Thus, $M_0 \cong U\oplus <4>$.

We consider its sublattice
$$
M_{II}=4M_0^\ast =\bz (4f_2) \oplus \bz f_3 \oplus \bz (4f_{-2})
\cong U(16)\oplus <4>.
$$
(here $S(n)$ denote a form $S$ which is multiplied by $n$).
We consider a subgroup $W^{(4)}(M_{II})$ generated by reflections in
$\Delta^{(4)}\subset M_{II}$ where $\Delta^{(4)}$ denote all elements with
the square $4$ (all these elements give
reflections of $M_{II}$). Then $W^{(4)}(M_{II})$
is a normal subgroup of finite
index in $O^+(M_0)$ with a fundamental polyhedron $\M_{II}$ with the
set of orthogonal vectors to $\M_{II}$ equals
$$
P(\M_{II})_{\pr}=\{\delta_1=-f_3,\ \delta_2=4f_2+f_3,\
\delta_3=4f_2+3f_3+4f_{-2},\ \delta_4=f_3+4f_{-2}\}
$$
where all elements $\delta_i$ have the square $4$ and the Gram matrix
$$
(\delta_i,\delta_j)=
\left(
\matrix
\hphantom{-}4 &  -4 & -12& -4 \\
 -4 & \hphantom{-}4 & -4 & -12\\
-12 &-4 & \hphantom{-} 4 & - 4 \\
 -4 & -12 & -4 & \hphantom{-}4
\endmatrix
\right).
\tag5.1
$$
Here $\M_{II}$ is a right quadrangle on the hyperbolic plane with all
its vertices at infinity.

Then
$$
A(P(\M_{II})_{\pr})\cong D_4
$$
is the group of symmetries of a right quadrangle and
$$
O^+(M_0)= W^{(4)}(M_{II}) \rtimes A(P(\M_{II})_{\pr}).
$$
We denote $P(\M_{II})=P(\M_{II})_{\pr}=\{\delta_1,\delta_2, \delta_3,
\delta_4\}$.

This group $W^{(4)}(M_{II})$ is also {\it elliptic} which
means (see \cite{N5}) that $\M_{II}$ has finite volume in the
hyperbolic plane $V^+(M_0)/\br_{++}$. Equivalently, for the cone
$\br _+Q_+=\br_+\delta_1+...+\br_+\delta_4$ and the corresponding dual
cone $\br_+Q^\ast_+=\{ x \in M_0\otimes \br\ |\
(x, \delta_i)\le 0\}$ we have similar to \thetag{2.2} sequence of embeddings
$$
\br_+ Q_+^\ast \subset \overline{V^+(M_0)} \subset \br_+Q_+.
\tag5.2
$$
It follows that $\br_+Q^\ast_+=\br_+\M_{II}$.
Moreover, $P(\M_{II})$ has a lattice Weyl vector $\rho$
which is defined by the property
$$
(\rho, \delta_i)=-(\delta_i,\delta_i)/2=-2.
$$
One can easily calculate that
$$
\rho =f_2+(1/2)f_3+f_{-2}.
\tag5.3
$$

We consider the corresponding ``complexified cone"
$$
\Omega(V^+(M_0))=M_0\otimes \br +iV^+(M_0)
$$
with coordinates $^t(z^\prime_1,z^\prime_2,z^\prime_3)$ where
$z^\prime=z^\prime_3f_2+z^\prime_2f_3+z^\prime_1f_{-2}
\in \Omega(V^+(M_0))$.
We consider another hyperbolic plane
 $\bz f_1 \oplus \bz f_{-1} \cong U$ with
$(f_1,f_1)=(f_{-1},f_{-1})$ $=0$, $(f_1,f_{-1})=-1$, and
its orthogonal sum $L$ with $M_0$. Thus, we can consider
$\Omega (V^+(M_0))$ as a cusp of the domain
$$
{\Cal H}_+ =
\{Z^\prime=^{t} \bigl((2 (z^\prime_2)^2-z^\prime_1z^\prime_3),\,
z^\prime_3,\,z^\prime_2,
\,z^\prime_1,\,1\bigr)\cdot z^\prime_0\in
\Cal H^{IV}\ |\ \hbox{Im}\,(z^\prime_1)>0\}.
$$
Below we  construct an automorphic form with a multiplier system
$F_2(Z^\prime)$ on ${\Cal H}_+$ with
respect to the group $O^+(L)$ which has the following properties.
By construction, this form will be anti-invariant relative to
$W^{(4)}(M_{II}) \rtimes A(P(\M_{II}))$
which  means that
$$
F_2(w.a(Z^\prime))=\det (w) F_2(Z^\prime)\ \text{for}\ \
w\in W^{(4)}(M_{II}),\ a \in A(P(\M_{II})).
\tag5.4
$$
By \thetag{5.2}---\thetag{5.4} and \thetag{5.11}, \thetag{5.12}
below, the form $F_2(Z^\prime)$ will have the Fourier expansion
with the properties:
\pagebreak
$$
\hskip-200pt
F_2(Z^\prime)=
$$
$$
\sum_{w\in W^{(2)}(M_{II})}{ \hskip-10pt \det(w)
\left(\exp(-\frac{\pi i}{2} (w(\rho),z^\prime))
- \hskip-10pt  \sum_{ a \in M_{II} \cap \br_{++}\M_{II}}
{\hskip-10pt m(a)\exp(-\frac{\pi i}{2} (w(\rho+a),z^\prime ))}\right)},
\tag5.5
$$
where $z^\prime = z^\prime_3f_2+z^\prime_2f_3+
z^\prime_1f_{-2} \in M_0 \otimes \br +i V^+(M_0)=
M_{II} \otimes \br +i V^+(M_{II})$ and
$m(a) \in \bz$.
For any primitive $a_0 \in M_{II}\cap \br_{++}\M_{II}$
with $(a_0,a_0)=0$, we will have the identity
$$
1- \sum_{t \in \bn}{m(ta_0)q^t}=
\prod_{k \in \bn}{(1-q^k)^{\tau(ka_0)}}
\tag5.6
$$
of power series of one variable $q$, where all $\tau(ka_0)\in \bz$.
More exactly,
$$
\tau(a)=3\ \ \text{for any}\ \  a \in M_{II}\cap \br_{++}\M_{II}\
\text{with}
\ (a,a)=0.
\tag5.7
$$

Like in \S~3, we use Fourier coefficients \thetag{5.5} and \thetag{5.6}
to construct the set of simple roots
$$
{}_s\Delta = {}_s\Delta^{\re}\cup {}_s\Delta^{\im}
$$
where ${}_s\Delta^{\re}=P(\M_{II})$ and
$$
\split
{}_s\Delta^{\im}&=
\{m(a)a\ |\ a\in M_{II}\cap \br_{++}\M_{II},\ (a,a)<0 \}\\
&\cup\{\tau(a)a\ |\ a\in M_{II}\cap \br_{++}\M_{II},\ (a,a)=0 \},
\endsplit
\tag5.8
$$
and the corresponding generalized Kac--Moody superalgebra
without odd real simple roots
$\geg^{\prime\prime}(M_{II}, {}_s\Delta)$. By construction,
the denominator function $\Phi (z^\prime)$ of
this algebra and the
automorphic form $F_2(Z^\prime)$ are connected by the formula
$$
\Phi (z^\prime )=F_2(4Z^\prime).
\tag5.9
$$
Below (see \thetag{5.14})
we construct the automorphic form $F_2(Z^\prime)$ and give its
development in the infinite product.
Like in \S~3 (see \thetag{3.2} and \thetag{3.3}), this gives
multiplicities for Weyl--Kac product formula of
the denominator function
$\Phi (z^\prime)$.

\vskip10pt

Now we work with the Siegel upper-half plane:
$\quad Z=\pmatrix z_1 &z_2\\z_2&z_3\endpmatrix \in \bh$.
Let
$$
\psi_{2,\frac1{2}}(z_1,z_2)=\eta(z_1)^3\vartheta_{11}(z_1,z_2).
$$
This  is a  Jacobi cusp  form of weight $2$ and index
$\dsize\frac 1{2}$ with respect to the  full Jacobi group and
with a multiplier system
$v_2:\Gamma^J\to \{\pm 1,\ \pm i\}$.
It means that the function
$\widetilde{\psi}_{2,\frac1{2}}(Z)=
\psi_{2,\frac1{2}}(z_1,z_2)\exp{(\pi i\,z_3)}$
satisfies the functional equations of types \thetag{4.4}--\thetag{4.5}
with $k=2$. We have the following exact formulae for the value
of multiplier system (see \thetag{4.2}, \thetag{4.3}):
$$
v_2(h(0,1,0))=v_2(h(1,0,0))=-1,\quad
v_2(\pmatrix 1&1\\0&1\endpmatrix^*)=
v_2(\pmatrix 0&-1\\1&0\endpmatrix^*)=i.
\tag5.10
$$
In accordance with  \thetag{4.6} and the well known
formula for $\eta(z_1)^3$, we have
$$
\align
{}&\psi_{2,\frac{1}2}(z_1, z_2)=q^{1/4}r^{-1/2}
\prod_{n\ge 1}\,(1-q^{n-1} r)(1-q^n r^{-1})(1-q^n)^4\\
{}&=\biggl(
\sum\Sb n\equiv 1\,mod\,2\endSb
\,(-1)^{\frac {n+1}2}
\exp{(\frac{\pi i\, }{4} n^2 z_1+\pi i\,nz_2)}
\biggr)
\times
\biggl(
\sum_{m\equiv 1\,mod\, 4} m\, \exp{(\frac{\pi i\,}{4} m^2 z_1)}
\biggr)
\\
{}&\hskip2cm=\sum
\Sb n\equiv 1\,mod\,4\\
\vspace{0.5\jot}
l\equiv 1\,mod\, 2\\
\vspace{0.5\jot}
2n>l^2
\endSb
c(n,l)\,\exp{(\frac{\pi i}{2}\,(nz_1+2lz_2))}.
\tag5.11
\endalign
$$

Analogues to \cite{G3} -- \cite{G4}, we can construct
the arithmetical lifting of
the Jacobi form
$\psi_{2,\frac{1}2}$ to a modular form of weight $2$
(with a multiplier system) with respect to the paramodular group
$$
\Gamma_2:=
\left\{\pmatrix *&2*&*&*\\
                    *&*&*&2^{-1}*\\
                    *&2*&*&*\\
                    2*&2*&2*&* \endpmatrix
\in Sp_4(\Bbb Q),\quad\text{all } *\in \bz\right\}.
$$
This group is  conjugate to the integral symplectic group
keeping
\newline
the skew-symmetric form
$$
\pmatrix 0&T_2\\-T_2&0\endpmatrix\qquad
T_2=\pmatrix 1&0\\0&2\endpmatrix.
$$
Using the lifting constructed in \cite{G3}, we obtain a cusp  form
$$
\multline
F_2(Z)=\text{Lift\,}\bigl(\widetilde\psi_{2,\frac{1}{2}}\bigr)(Z)\\
\hskip-20pt=\sum\Sb  n,m\equiv 1\, mod\, 4\\
\vspace{0.5\jot}
l\equiv 1\, mod\,  2\\
\vspace{0.5\jot}
 n,m>0,\ 2mn-l^2>0
\endSb
\sum_{d\,|\,(n,l,m)} d \biggl(\frac {-4}{d}\biggr)
\, c(\frac{nm}{d^2},\,\frac{l}{d})
\exp{(\frac{\pi i}{2}\,( nz_1+2lz_2+2mz_3))},
\endmultline
\tag5.12
$$
where
$$
\biggl(\frac {-4}{d}\biggr)
=\cases 1&\text{if }d\equiv 1\mod 4\\
-1&\text{if }d\equiv -1\mod 4.
\endcases
$$
The multiplier system of $F_2$ is defined by \thetag{5.10} and
the equation
$$
v_2(\pmatrix {}^tU&0\\0&U\endpmatrix)=1
\qquad\text{with }\  U=\pmatrix 0&{\sqrt 2}^{-1}\\ \sqrt 2&0\endpmatrix.
$$
To define an expansion in the infinite product we introduce
a weak Jacobi form of weight zero and index 2
$$
\align
\phi_{0,2}(z_1,z_2)&=\frac 1{288\Delta_{12}(z_1)}
\bigl(E_4 (z_1)E_{4,1}^2(z_1,z_2)-E_{6,1}^2(z_1,z_2)\bigr)\\
\vspace{2\jot}
{}&=\sum\Sb n,l \endSb\  f_2(n,l)\exp{(2\pi i\, n z_1 +l z_2)}.
\tag5.13
\endalign
$$
Combining the method of \cite{G3} with
the method used above, we get
$$
F_2(Z)=\exp{(\frac{\pi i}2\,(z_1-2z_2+2z_3))}\kern-11pt
\prod\Sb n,\,l,\,m \in \Bbb Z\\ (n,l,m)>0\endSb
\bigl(1-\exp{(2\pi i\,(nz_1+lz_2+2mz_3))}\bigr)^{f_2(nm,l)},
$$
where $(n,l,m)>0$ means that $n\ge 0$, $m\ge 0$,
$l$ is an  arbitrary integral if $n+m>0$, and $l>0$ if $n=m=0$.

\vskip10pt

For any  paramodular symplectic group the analogue of
Lemma 1.1 is true (see  \cite{G2}). Using this isomorphism,
one can easily rewrite the result about $F_2$ above
and obtain the analogue of the formula \thetag{4.13}
in terms of the orthogonal group:
$$
\multline
F_2(Z^\prime)=\sum
\Sb n,m\equiv 1\, mod\,4\\
l\equiv 1\, mod\,2
\endSb\
\sum_{d|(n,l,m)} d \biggl(\frac {-4}{d}\biggr)\,
c(\frac{nm}{d^2},\,\frac{l}{d})
\exp{(\frac{\pi i}{2}\,( nz'_1+2lz'_2+mz'_3))}=\\
\vspace{2\jot}
\exp{(\frac{\pi i}2\,(z'_1-2z'_2+z'_3))}\,
\prod\Sb n,\,l,\,m\in \Bbb Z\\ (n,l,m)>0\endSb
\bigl(1-\exp{(2\pi i\,(nz'_1+lz'_2+mz'_3))}\bigr)^{f_2(nm,l)},
\endmultline
\tag5.14
$$
where $c(\cdot,\cdot)$ and $f_2(\cdot,\cdot)$ are
the Fourier coefficients of the corresponding
Jacobi forms of weight $2$ and $0$
and $z'_i$ are variables from the homogeneous domain of the
orthogonal group. This gives the formula \thetag{5.9}.

\head
\S~6. Appendix: Generalized Kac--Moody superalgebras without
odd real simple roots
\endhead

We refer to the paper of V. Kac \cite{Ka2} for basic definitions related with
Lie superalgebras.

We fix a complex $n\times n$
matrix $A=(a_{ij})$ where $i,j \in I=\{1,2,...,n\}$.
Moreover, we fix a subset $\tau \subset I$.

Like in  \cite{Ka1}, one can define a graded Lie superalgebra
$\tilde\geg (A, \tau)$, its maximal graded ideal $\rr$
trivially intersecting the Cartan subalgebra $\hh$ and Lie superalgebra
$\geg(A, \tau)=\tilde\geg(A, \tau)/\rr$. The only difference is
(see \cite{Ka1, Ch. 1}) that
the generators $\hh$ and generators $e_i$, $f_i$, $i \in I-\tau$,
are even,
and generators $e_i$, $f_i$, $i \in \tau$, are odd.
For $\tau=\emptyset$, one obtains Chevalley Lie algebras $\geg (A)$ which
were considered in \cite{Ka1}.

The matrix $A$ above is called a {\it generalized generalized
Cartan matrix} if $A$ is a real matrix which satisfies the
following conditions:

\vskip5pt

(C1$^\prime$) either $a_{ii}=2$ or $a_{ii}\le 0$;

(C2$^\prime$) $a_{ij}\le 0$ if $i\not=j$, and $a_{ij}\in \bz$
if $a_{ii}=2$;

(C3$^\prime$) $a_{ij}=0$ implies $a_{ji}=0$.

\vskip5pt

If $A$ is a generalized generalized Cartan matrix, then
the Lie superalgebra $\geg(A, \tau)$ is called a
{\it generalized Kac--Moody superalgebra}. We assume
the additional condition

\vskip5pt

(C4$^\prime$) if $i\in \tau$, then $a_{ii}\le 0$.

\vskip5pt

If $A$ satisfies (C4$^\prime$), then
the generalized Kac--Moody superalgebra $\geg(A, \tau)$
is called a {\it generalized Kac--Moody superalgebra without odd real
simple roots}.
In particular, for $\tau=\emptyset$ one obtains a
generalized Kac--Moody algebra $\geg (A)$.

Results of the classical book by V. Kac \cite{Ka1}
are divided in 3 types:

\vskip5pt

{\it Type 1.}
Results which are valid for Lie algebras
$\geg (A)$ with
an arbitrary or an arbitrary symmetrizable
matrix $A$. These are results of Ch. 1 (Basic definitions); Ch. 2
(The invariant bilinear form and the generalized Casimir operator);
Ch. 9 (Highest weight modules over Lie algebra $\geg(A)$).

{\it Type 2.}
Results which are valid for an arbitrary or an arbitrary
symmetrizable (i.e. $A$ is symmetrizable) generalized Kac--Moody
algebra $\geg(A)$. These are results of
Ch. 3 (Integrable representations and the Weyl group of a Kac--Moody
algebra); Ch. 4 (Some properties of generalized Cartan matrices);
Ch. 5 (Real and imaginary roots); Ch. 10 (Integrable highest weight
modules: the character formula);
Ch. 11 (Integrable highest weight modules:
the weight system, the contravariant Hermitian form and the
restriction problem). We mention that \cite{Ka1, \S~11.13}
contains changes one should make to adapt results of these Chapters (where
Kac--Moody algebras (i.e. $a_{ii}=2$ for any $i \in I$)
are considered) to generalized Kac--Moody
algebras. See also \cite{Bo1} and \cite{Bo2} where generalized
Kac--Moody algebras were first introduced and studied,
but we prefer to follow to \cite{Ka1}.

{\it Type 3.} Results on affine Kac--Moody algebras.
Chs. 6, 7, 8, 12, 13, 14.

\vskip5pt

The general remark is that all type 1 results may be
adapted to Lie superalgebras
$\geg(A,\tau)$ with arbitrary or symmetrizable
complex matrix $A$.
All type 2 results may be adapted to Kac--Moody superalgebras
without odd real simple roots $\geg (A,\tau)$.
All type 3 results may be adapted to appropriate
affine generalized Kac--Moody superalgebras without
odd real simple roots $\geg (A,\tau)$.
Here we consider changes in \cite{Ka1} one should make for
these adaptations of Type 1 and Type 2 results.

\vskip5pt

Type 1 results: One should introduce the following changes in
Chs. 1, 2 and 9 of \cite{Ka1} to adapt them to
Lie superalgebras $\geg(A, \tau)$ with arbitrary complex $A$:

Instead of \cite{Ka1, (1.3.3)}, one has:
$$
\geg_{\alpha_i}=\bc e_i,\
\geg_{-\alpha_i}=\bc f_i,\
\geg_{s\alpha_i}=0,\  \text{if}\  |s|>1
\tag6.1
$$
for $i \in I-\tau$ and for $i \in \tau$ and $a_{ii}=0$;
$$\gather
\geg_{\alpha_i}=\bc e_i,\ \ \geg_{-\alpha_i}=\bc f_i,\ \ \geg_{s\alpha_i}=0,
 \text{ if }\  |s|>2,\\
\geg_{2\alpha_i} = \bc [e_i,e_i]\not=\{ 0\}, \ \
\geg_{-2\alpha_i}=\bc [f_i,f_i]\not=\{0\}
 \tag6.2
\endgather
$$
for $i \in \tau$ and $a_{ii}\not=0$.

Instead of \cite{Ka1, Lemma 1.3}, one has: If
$\beta\in \Delta_+\ -\{\alpha_i,
2\alpha_i\}$, then $(\beta+\bz\alpha_i)\cap \Delta\subset \Delta_+$.

Instead of Chevalley involution $\omega$ (and the corresponding
involution $\tilde\omega$) one should consider the
Chevalley automorphism of the period $4$ which is defined by
the property:
$e_i\mapsto -(-1)^{\deg e_i}f_i$, $f_i\mapsto -e_i$\
($i \in I$), $h\mapsto -h$\  ($h \in \hh$).
It has $\omega^2$ identical on
$\geg(A,\tau)_\0o$ and equal $-1$ on $\geg(A,\tau)_\1o$.

There are the following changes in \cite{Ka1, Ch. 9}.
Let $V$ be a $\geg(A,\tau)$-module from the category $\Cal O$ and
$V=\bigoplus_{\lambda \in \hh^\ast}{V_\lambda}$.
There are two possibile definitions of the formal character of $V$:
$$
\ch V=\sum_{\lambda \in \hh^\ast}{(\dim V_\lambda)e(\lambda)};
\tag6.3
$$
and
$$
\cha V=\sum_{\lambda \in \hh^\ast}
{(\dim V_{\lambda,\0o}-\dim V_{\lambda,\1o} )e(\lambda)}.
\tag6.4
$$
For a Verma module $M_{\overline{i}}(\Lambda)$ with the highest weight
vector of the degree $\overline{i}\in \{\overline0,\overline1\}$ and
$\Lambda \in \hh^\ast$, one has
$$
\multline
\ch M_{\overline{i}}(\Lambda )=\\
e(\Lambda )\prod_{\alpha \in \Delta_{+,\0o}}{(1-e(-\alpha))^{-\mult~\alpha}}
\prod_{\alpha \in \Delta_{+,\1o}}{(1+e(-\alpha))^{-\mult~\alpha}}
= e(\Lambda)\ch M_\0o (0)
\endmultline
\tag6.5
$$
and
$$
\cha M_{\overline{i}}(\Lambda )=
(-1)^{\overline{i}}e(\Lambda )\prod_{\alpha \in \Delta_+}
{(1-e(-\alpha))^{-\mult~\alpha}}
=(-1)^\io e(\Lambda)\cha M_\0o (0)
\tag6.6
$$
where for both formulae \thetag{6.5} and \thetag{6.6}
$$
\mult~\alpha =
\cases
\ \dim \geg(A,\tau )_\alpha &\text{if $\alpha \in \Delta_{\0o}$}, \\
-\dim \geg(A,\tau )_\alpha &\text{if $\alpha \in \Delta_{\1o}$}.
\endcases
\tag6.7
$$
One should consider only homomorphisms of $\geg(A, \tau)$-modules
which preserve even and odd parts of modules.

\vskip5pt

Type 2 results:
In \cite{Ka1, \S~11.13} one can find changes one should make to
adapt Type 2 results of \cite{Ka1} to generalized Kac--Moody
algebras (i. e. algebras $\geg (A)$ where $A$ satisfies the
conditions (C1$^\prime$),
(C2$^\prime$), (C3$^\prime$) (and $\tau=\emptyset$). We give
here additional changes one should make to adapt these results
to generalized Kac--Moody superalgebras
without odd real simple roots $\geg (A,\tau)$
(i.e. with condition (C4$^\prime$) instead of $\tau=\emptyset$).

Instead of \cite{Ka1, Lemma 11.13.1}, one has:

\proclaim{Statement 6.1}
Let $V$ be an integrable $\geg (A)$-module,
$\lambda \in \hh^\ast$ and $0\not= v \in V_\lambda$ such that
$f_i(v)=0$. Then for $k\ge 1$,
$$
f_ie_i^k(v)=
\cases
-k(\langle \lambda, \alpha_i^\vee\rangle +
\frac{1}{2}(k-1)a_{ii})e_i^{k-1}(v)
&\text{if $i\in I-\tau$;} \\
-[k/2]a_{ii}e_i^{k-1}(v),
&\text{if $i\in \tau$ and $k \equiv 0 \mod 2$;}\\
-([k/2]a_{ii}+\langle\lambda, \alpha_i^\vee\rangle)e_i^{k-1}(v),
&\text{if $i\in \tau$ and $k \equiv 1 \mod 2$.}\\
\endcases
$$
In particular, $e_i^k(v)\not=0$ for all $k=1,2,...$ if
$a_{ii}\le0$, $\langle\lambda, \alpha_i^\vee\rangle<0$ and either
$i \in I-\tau$, or $i \in \tau$ and $a_{ii}<0$.
\endproclaim

Using this statement, one obtains the following statement
instead of \cite{Ka1, Corollary 11.13.1}:

\proclaim{Statement 6.2}
Assume that $a_{ii} \le 0$ and $i \in I-\tau$ for $a_{ii}=0$.
Assume that $\text{supp}(\alpha+\alpha_i)$ is connected. Then
$\alpha+k\alpha_i\in \Delta_+$ for all $k \in \bz_+$.
\endproclaim

All results of \cite{Ka1, Ch. 4} are valid for generalized
Kac--Moody superalgebras except that there exists one
additional zero $1\times1$ matrix which gives a generalized
Kac--Moody algebra for $\tau=\emptyset$ and a generalized
Kac--Moody superalgebra for $\tau=\{1\}$.

To describe like in \cite{Ka1, Ch. 5} the set of roots $\Delta$, we set
$\tau_{<0}=\{i \in \tau\ |\ a_{ii}<0\}$,
$\tau_0=\{i \in \tau\ |\ a_{ii}=0\}$, and
$$
\Pi^{\im}_\0o=\{\alpha_i \in \Pi^{\im}\ |\ i \in I-\tau \} ;
\hskip20pt
\Pi^{\im}_\1o=\{\alpha_i \in \Pi^{\im}\ |\ i \in \tau \} ;
$$
$$
\Pi^{\im}_{\1o,<0}=\{\alpha_i \in \Pi\  |\ i \in \tau_{<0} \} ;
\hskip20pt
\Pi^{\im}_{\1o,0}=\{\alpha_i \in \Pi\ |\ i \in \tau_0 \}.
$$
Let
$$
\split
K=&\{\alpha \in Q_+\ |\ \alpha \in -C^\vee\ \text{and supp $\alpha$ is
connected} \}\\
-&\left((\bigcup_{k \ge 2}{k \Pi^{\im}_\0o})\cup
(\bigcup_{k \ge 3}{k\Pi^{\im}_{\1o,<0}})\cup
(\bigcup_{k \ge 2}{k\Pi^{\im}_{\1o,0}})\right).
\endsplit
$$
Then we have
$$
\Delta_+^{\im}\subset \bigcup_{w\in W}{w(K)}.
$$
Like in \cite{Ka1, (11.13.3)} and using Statement 6.2,
we get the opposite inclusion and thus

\proclaim{Statement 6.3}
$$
\Delta_+^{\im} = \bigcup_{w\in W}{w(K)}.
$$
if $\Pi^{\im}_{\1o, 0}=\emptyset$ (i.e. $\tau_0=\emptyset$).
\endproclaim

We don't know what will be the analog of Statement 6.3
in general when $\Pi^{\im}_{\1o, 0}$ is not empty.

Using \cite{Ka1, Lemma 11.13.2} and superalgebras analog of
\cite{Ka1, Proposition 9.11}, it follows that
symmetrizable (i.e. with symmetrizable $A$)
generalized Kac--Moody superalgebras without
odd real simple roots $\geg (A,\tau)$
have similar defining relations as generalized
Kac--Moody algebras $\geg(A)$.

\proclaim{Statement 6.4} A symmetrizable generalized Kac--Moody
superalgebra without odd real simple roots $\geg (A,\tau)$
is a Lie superalgera with even generators
$\hh$, $e_i$, $f_i$, $i \in I-\tau$, and odd generators
$e_i$, $f_i$, $i \in \tau$, and defining relations
$$
\split
&[e_i,f_j]=\delta_{ij}\alpha_i^\vee \ \ (i,j \in I), \\
&[h,h^\prime]=0 \ \ (h, h^\prime \in \hh),\\
&[h,e_i]=\langle\alpha_i,h\rangle e_i, \\
&[h, f_i]=-\langle\alpha_i,h\rangle f_i,\ \ (i \in I;\ h \in \hh),
\endsplit
$$
$$
(\text{ad\ }e_i)^{1-a_{ij}}e_j=0,\ (\text{ad\ }f_i)^{1-a_{ij}}f_j=0,\
\text{if}\ a_{ii}=2\ \text{and}\ i \not=j,
$$
$$
[e_i,e_j]=0,\ [f_i,f_j]=0,\ \text{if}\ a_{ij}=0.
$$
Remark that the last relation is not trivial also for $i=j$ if
$a_{ii}=0$ and $i \in \tau$.
\endproclaim

Let us consider a diagonal generalized generalized Cartan matrix
$B=(b_{ij})$ with $b_{ii}\le 0$ for any $i \in I$ and some
$\tau \subset I$, and the generalized
Kac--Moody superalgebra $\geg (B, \tau)$.
Using \thetag{6.1}, \thetag{6.2}, one obtains for $\geg (B, \tau)$:
$$
\multline
\ch L_\0o (0)=(\ch M_\0o (0))^{-1}\\
=\prod_{i \in I-\tau}
{\hskip-5pt \left(1-e(-\alpha_i)\right)}
\prod_{i \in \tau_{<0}}
{\hskip-5pt \left(\left(1-e(-2\alpha_i)\right)
\left(1+e(-\alpha_i)\right)^{-1}\right)}
\prod_{i \in \tau_0}
{\hskip-5pt \left(1+e(-\alpha_i)\right)^{-1}}\\
= \prod_{i \in I-\tau_0}{\left(1-e(-\alpha_i)\right)}
\prod_{i \in \tau_0}{\left(\sum_{k\ge 0}{(-1)^k e(-k\alpha_i)}\right)}.
\endmultline
\tag6.8
$$
It follows that  $\ch L_\0o (0)$ is the sum
$$
\ch L_\0o (0) = \sum{\epsilon (s)e(-s)},
\tag6.9
$$
over all sums of simple roots $s$. Here
the sign $\epsilon (s)=(-1)^n$ if $s$ is a sum of
$n$ simple pairwise perpendicular imaginary simple roots $\alpha_i$
which are distinct if $i \in I-\tau_0$; and $\epsilon (s)=0$ otherwise.

Similarly, for the generalized generalized Cartan matrix $B$
above and the character $\cha$, one has
$$
\multline
\cha L_\0o (0)=(\cha M_\0o(0))^{-1}\\
=\prod_{i \in I-\tau}{\hskip-5pt \left(1-e(-\alpha_i)\right)}
\prod_{i \in \tau_{<0}}{\hskip-5pt
\left(\left(1-e(-2\alpha_i)\right)\left(1-e(-\alpha_i)\right)^{-1}\right)}
\prod_{i \in \tau_0}{\hskip-5pt \left(1-e(-\alpha_i)\right)^{-1}}\\
= \prod_{i \in I-\tau}{\left(1-e(-\alpha_i)\right)}
\prod_{i \in \tau_{<0}}{\left(1+e(-\alpha_i)\right)}
\prod_{i \in \tau_0}{\left(\sum_{k\ge 0}{e(-k\alpha_i)}\right)} .
\endmultline
$$
It follows that  $\cha L_\0o (0)$ is the sum
$$
\cha L_\0o(0) = \sum{\epsilon (s)e(-s)},
\tag6.11
$$
over all sums of simple roots $s$. Here
the sign $\epsilon (s)=(-1)^{n_\0o}$ if $s$ is a sum of
pairwise perpendicular
$n_\0o$ imaginary simple roots $\alpha_i$, $i \in I-\tau$,
and $n_\1o$ imaginary simple roots $\alpha_i$, $i \in \tau$,
which are distinct if $i \in I-\tau_0$; and $\epsilon (s)=0$ otherwise.

Now we consider an arbitrary generalized generalized
Cartan matrix $A$.
We consider an irreducible $\geg (A, \tau)$-module
$L_\io (\Lambda )$ with a heighest weight vector of
the degree $\io\in \{\0o,\1o\}$ and $\Lambda\in \hh^\ast$
which satisfies the condition
$$
\langle\Lambda, \alpha_i^\vee\rangle \in \bz_+\ \text{if $a_{ii}=2$},\
\langle\Lambda, \alpha_i^\vee\rangle \ge 0 \  \text{for all $i \in I$}.
\tag6.12
$$

We remind that $\rho \in \hh^\ast$ is defined by relations
$$
\langle\rho, \alpha_i^\vee\rangle=\frac{1}{2}a_{ii},\ \
\text{for all\ } i \in I.
$$
Equivalently, for a symmetrizable $A$,
$$
(\rho|\alpha_i)=\frac{1}{2}(\alpha_i|\alpha_i)\ \
\text{for all\ } i \in I.
$$

 From \thetag{6.5},
\thetag{6.6}, \thetag{6.8} --- \thetag{6.11},
like in \cite{Ka1, \S~11.13}, we get
Weyl--Kac--Borcherds  character formulae for symmetrizable
generalized Kac--Moody superalgebras without
odd real simple roots $\geg (A,\tau)$.
Formulae \thetag{6.8} --- \thetag{6.11} are particular cases of them.

\proclaim{Statement 6.5}
$$
\ch L_\io (\Lambda )=
\sum_{w \in W}{\epsilon (w)w(S_{\Lambda,\io})}\Big/
e(\rho )\prod_{\alpha \in \Delta_{+,\0o}}
{\hskip-5pt (1-e(-\alpha))^{\mult~\alpha}}
\prod_{\alpha \in \Delta_{+, \1o}}
{\hskip-5pt (1+e(-\alpha))^{\mult~\alpha }}
\tag6.13
$$
where $S_{\Lambda, \io}$ is the sum
$$
S_{\Lambda,\io} = e(\Lambda+\rho)\sum{\epsilon (s)e(-s)}
\tag6.14
$$
over all sums of imaginary simple roots $s$. Here
the sign $\epsilon (s)=(-1)^n$ if $s$ is a sum of
$n$ pairwise perpendicular and perpendicular to $\Lambda$
imaginary simple roots $\alpha_i$
which are distinct for $i \in I-\tau_0$; and $\epsilon (s)=0$ otherwise.
\endproclaim

\proclaim{Statement 6.6}
$$
\cha L_\io (\Lambda )=
\sum_{w\in W}
{\hskip-5pt \epsilon (w)w(S_{\Lambda,\io})} \Big/
e(\rho) \hskip-5pt \prod_{\alpha \in \Delta_+}
{\hskip-5pt (1-e(-\alpha ))^{\mult~\alpha}}
\tag6.15
$$
where $S_{\Lambda,\io}$ is the sum
$$
S_{\Lambda,\io} = (-1)^{\io}e(\Lambda+\rho)\sum{\epsilon (s)e(-s)}
\tag6.16
$$
over all sums of imaginary simple roots $s$. Here
$\epsilon (s)=(-1)^{n_\0o}$ if $s$ is a sum of
$n_\0o+n_\1o$ pairwise perpendicular and perpendicular
to $\Lambda$  $n_\0o$ imaginary simple roots $\alpha_i$,
$i \in I-\tau$, and $n_\1o$ imaginary simple roots $\alpha_i$,
$i \in \tau$, which are distinct for $i \in I-\tau_0$;
and $\epsilon (s)=0$ otherwise.
\endproclaim

In particular, for $\Lambda=0$ and $\io=\0o$ we get the denominator
identities for $\geg (A, \tau)$.

\proclaim{Statement 6.7}
$$
e(\rho )\prod_{\alpha \in \Delta_{+,\0o}}
{\hskip-5pt (1-e(-\alpha))^{\mult~\alpha}}
\prod_{\alpha \in \Delta_{+, \1o}}
{\hskip-5pt (1+e(-\alpha))^{\mult~\alpha }}=
\sum_{w \in W}{\epsilon (w)w(S)}
\tag6.17
$$
where $S$ is the sum
$$
S = e(\rho)\sum{\epsilon (s)e(-s)}
\tag6.18
$$
over all sums of imaginary simple roots $s$. Here
the sign $\epsilon (s)=(-1)^n$ if $s$ is a sum of
$n$ pairwise perpendicular
imaginary simple roots $\alpha_i$
which are distinct for $i \in I-\tau_0$; and $\epsilon (s)=0$ otherwise.
\endproclaim

\proclaim{Statement 6.8}
$$
e(\rho) \hskip-5pt \prod_{\alpha \in \Delta_+}
{\hskip-5pt (1-e(-\alpha ))^{\mult~\alpha}}
=
\sum_{w\in W}{\hskip-5pt \epsilon (w)w(S)}
\tag6.19
$$
where $S$ is the sum
$$
S = e(\rho)\sum{\epsilon (s)e(-s)}
\tag6.20
$$
over all sums of imaginary simple roots $s$. Here
$\epsilon (s)=(-1)^{n_\0o}$ if $s$ is a sum of
pairwise perpendicular
$n_\0o$ imaginary simple roots $\alpha_i$, $i \in I-\tau$, and $n_\1o$
imaginary simple roots $\alpha_i$, $i \in \tau$, which are distinct
for $i \in I-\tau_0$; and $\epsilon (s)=0$ otherwise.
\endproclaim


Instead of a generalized Kac--Moody superalgebra
without odd real simple roots $\geg (A,\tau)$
one can consider the derived algebra
$\geg^\prime (A, \tau) = [\geg (A, \tau), \geg (A, \tau)]$ which is
also called generalized Kac--Moody superalgebra without
odd real simple roots.
By Statement 6.4, a symmetrizable $\geg ^\prime (A, \tau)$ is
generated by $h_i=\alpha_i^\vee$, $e_i$, $f_i$, $i \in I$,
where all $h_i$, $i \in I$, are even, $e_i$, $f_i$ are even if
$i \in I-\tau$, and $e_i$, $f_i$ are odd if $i \in \tau$, and has
defining relations ($i, j \in I$)
$$
[e_i,f_j]=\delta_{ij}h_i,\  [h_i,h_j]=0 \ [h_i,e_j]=a_{ij}e_j, \
[h_i, f_j]=-a_{ij}f_j;
$$
$$
(\text{ad\ }e_i)^{1-a_{ij}}e_j=0,\  (\text{ad\ }f_i)^{1-a_{ij}}f_j=0\ \
\text{if}\ a_{ii}=2\ \text{and}\ i \not=j;
$$
$$
[e_i,e_j]=0,\ [f_i,f_j]=0 \ \ \text{if}\ a_{ij}=0.
\tag6.21
$$
We remind (see \cite{Ka1}) that
this algebra is graded by the root lattice
$$
Q=\bigoplus_{i \in I}{\bz \alpha_i},
\tag6.22
$$
where
$$
\deg e_i=\alpha_i=-\deg f_i, \deg h_i=0.
\tag6.23
$$
There is a pairing
between $Q$ and the coroot lattice
$$
Q^\vee=\bigoplus_{i \in I}{\bz h_i}
\tag6.24
$$
which is defined by
$$
\langle h_i, \alpha_j \rangle=a_{ij}.
\tag6.25
$$
Remind that $\hh^\prime = Q^\vee \otimes \bc$ is Cartan subalgebra of
$\geg^\prime (A, \tau)$.
Like in \cite{Ka1}, all results above (Statements 6.1 --- 6.3, 6.5---6.8
and \thetag{6.21}) are valid for
$\geg^\prime (A, \tau)$. Advantage of considering $\geg^\prime (A, \tau)$
is that one can consider also infinite matrices $A$ with a countable set of
indexes $I$. One can consider $\geg^\prime (A, \tau)$
with a countable $I$ as a union of
$\geg ^\prime (A^\prime, \tau^\prime )$ with finite $I^\prime$,
$I^\prime \subset I$ and $\tau^\prime=\tau\cap I^\prime$.

\vskip5pt

R. Borcherds \cite{Bo1} (see also \cite{Bo2})
uses another notation for generalized Kac--Moody algebras. One can use
similar notation to define generalized Kac--Moody superalgebras
without odd real simple roots as follows.

Suppose that we have (compare with \cite{Bo1})

\vskip5pt

(i) A real vector space $H$ with a symmetric bilinear inner
product ( , ).

(ii) A sequence of elements $h_i \in H$
indexed by a countable set $I$, such
that $(h_i, h_j) \le 0$ if $i\not=j$ and $2(h_i, h_j)/(h_i,h_i)$ is an
integer if $(h_i, h_i)>0$.

(iii) A subset $\tau \subset I$ such that $(h_i, h_i)\le 0$ for any
$i \in \tau$.
\vskip5pt

The matrix $(h_i, h_j)$, $i,j \in I$, is called  a
{\it symmetrized generalized generalized Cartan matrix}.
The generalized Kac--Moody superalgebra
$G$ associated to (i)---(iii) is defined
to be the Lie superalgebra generated by $H$ and elements
$e_i$, $f_i$, $i \in I$, where $H$ and $e_i$, $f_i$,
$i \in I-\tau$, are even, and $e_i$, $f_i$, $i \in \tau$, are odd. It has
the following defining relations:

\vskip5pt

(1) The image of $H$ in $G$ is commutative. (In fact the natural map from
$H$ to $G$ is injective so we can consider $H$ to be
an abelian subalgebra of G.)

(2) If $h$ is in $H$, then $[h,e_i]=(h, h_i)e_i$ and
$[h, f_i]=-(h, h_i)f_i$.

(3) $[e_i, f_i]=h_i$ if $i=j$, $0$ if $i \not=j$.

(4) If $(h_i, h_i)>0$ and $i \not=j$, then
$(\text{ad\ }e_i)^{1-2(h_i, h_j)/(h_i,h_i)}e_j=0$ and
\newline
$(\text{ad\ }f_i)^{1-2(h_i, h_j)/(h_i,h_i)}f_j=0$;

(5) If $(h_i, h_j)=0$, then $[e_i, e_j]=[f_i,f_j]=0$.
\vskip5pt

The algebra $G$ is graded by
$$
\widetilde{Q}=\bigoplus_{i \in I}{\bz r_i},
\tag6.26
$$
where
$$
\deg e_i=r_i=-\deg f_i, \ \ \deg h_i=0.
\tag6.27
$$

We set for $i, j\in I$,
$$
a_{ij}=
\cases
(h_i,h_j) &\text{if $(h_i, h_i)\le 0$},\\
2(h_i,h_j)/(h_i, h_i) &\text{if $(h_i, h_i)>0$}.
\endcases
$$
The matrix $A=(a_{ij})$ is a
symmetrizable generalized generalized Cartan matrix. It
defines the symmetrizable generalized Kac--Moody
superalgebra without odd real simple roots $\geg^\prime (A, \tau)$.
Below we denote its standard generators $h_i$,$e_i$, $f_i$
in \thetag{6.21} as $\widetilde{h}_i$, $\widetilde{e}_i$,
$\widetilde{f}_i$ respectively.

We have the homomorphism
$$
\pi: \geg^\prime (A, \tau) \to G,
\tag6.28
$$
which is defined by
$$
\pi (\widetilde{h}_i)=
\cases
h_i, &\text{if $(h_i,h_i)\le 0$},\\
(2/(h_i,h_i))h_i &\text{if $(h_i,h_i)>0$};
\endcases
$$
$$
\pi (\widetilde{e}_i)=
\cases
e_i, &\text{if $(h_i,h_i)\le 0$},\\
+\sqrt{2/(h_i,h_i)}e_i &\text{if $(h_i,h_i)>0$};
\endcases
$$
$$
\pi (\widetilde{f}_i)=
\cases
f_i, &\text{if $(h_i,h_i)\le 0$},\\
+\sqrt{2/(h_i,h_i)}f_i &\text{if $(h_i,h_i)>0$};
\endcases
$$
which is evidently homogeneouse for gradings
\thetag{6.22}, \thetag{6.23} and \thetag{6.28}, \thetag{6.29}.
This homomorphism
defines the isomorphism of $G$ with the
quotient of $\geg^\prime (A, \tau)$ by the ideal in the center
$\cc = \{ h \in \hh^\prime \ |\
(\pi (h), h_i)=0,\ \text{for any\ } i \in I\}$.
Using this isomorphism, one can transfer all results above
(Statements 6.1---6.3, 6.5---6.8 and \thetag{6.21}) to $G$.

\vskip5pt

Let us consider the canonical homomorphism of abelian groups
$$
\pi : \widetilde{Q}\to H\ \ (r_i \mapsto h_i).
\tag6.29
$$
If $\pi$ has finite preimages on the semigroup
$\widetilde{Q}_+\subset \widetilde{Q}$, where
$$
\widetilde{Q}_+=\bigoplus_{i \in I}{\bz_+r_i},
$$
and $\pi^{-1}(0)\cap \widetilde{Q}_+=0$,
we can replace the $\widetilde{Q}$ grading above by
the grading using
$$
Q=\pi (\widetilde{Q}) \subset H,\ \  Q_+=\pi (\widetilde{Q}_+)
\tag6.30
$$
letting for $\alpha\in Q$
$$
G_\alpha=\bigoplus_
 {\{ \widetilde{\alpha} \in \widetilde{Q} |
\pi(\widetilde{\alpha})=\alpha \}}
{G_{\widetilde{\alpha}}}.
$$
For this grading, $G_0=H$ and
$[h, x_\alpha]=(h, \alpha)x_\alpha$ for any $h \in H$,
$x_\alpha \in G_\alpha$, $\alpha \in Q$.
Here $G_\alpha=G_{\alpha,\0o}\oplus G_{\alpha,\1o}$ where
even and odd parts may be both non-zero. For $\alpha \in Q$ we set
$$
\mult_\io \alpha=(-1)^{\io}\dim G_{\alpha,\io},\ \
\mult~\alpha =\mult_\0o \alpha+\mult_\1o \alpha =
\dim G_{\alpha,\0o}-\dim G_{\alpha,\1o}.
\tag6.31
$$
An $\alpha \in Q$ is called a root if
$\dim G_\alpha =\dim G_{\alpha,\0o} + \dim G_{\alpha,\1o} >0$.
Let $\Delta\subset Q$ be the set of all roots,
$\Delta_\pm=\Delta \cap \pm Q_+$.

\vskip5pt

We formulate Statements 6.5---6.8 for $G$ with this grading.

Let $\Lambda \in H$ is such that $(\Lambda ,h_i) \ge 0$
and $(\Lambda ,h_i) \in \bz$ if $(h_i, h_i)>0$ for any $i \in I$.
Let $L(\Lambda)$ be an irreducible $G$-module with the highest
weight $\Lambda \in H \subset H^\ast$ and a highest weight vector of
even degree. For $\alpha \in Q$ we set
$L(\Lambda)_{\Lambda+\alpha}= U(G)_\alpha (L(\Lambda)_\Lambda)$
where $\dim L(\Lambda)_\Lambda=1$. See \cite{Ka1, \S~9.10}.
We remark that the $Q$-grading of $G$ above induces
the $Q$-grading of the universal enveloping algebra $U(G)$.
We set
$$
\ch L(\Lambda )=\sum_{\alpha \in Q}
{\dim L(\Lambda )_{\Lambda+\alpha }e(\Lambda+\alpha)}.
\tag6.32
$$
Similarly,
$$
\cha L(\Lambda )=\sum_{\alpha \in Q}
{(\dim L(\Lambda )_{\Lambda + \alpha ,\0o }
-\dim L(\Lambda )_{\Lambda +\alpha , \1o}) e(\Lambda+\alpha)}.
\tag6.33
$$

\proclaim{Statement 6.5$^\prime$}
$$
\ch L(\Lambda )=
\sum_{w \in W}{\epsilon (w)w(S_\Lambda)}\Big/
e(\rho )\prod_{\alpha \in \Delta_+}
{(1-e(-\alpha))^{\mult_\0o \alpha}}
\prod_{\alpha \in \Delta_+ }
{(1+e(-\alpha))^{\mult_\1o \alpha }}
$$
where $S_\Lambda$ is the sum
$$
S_\Lambda = e(\Lambda+\rho)\sum{\epsilon (s)e(-\pi (s))}
$$
over all sums of imaginary simple roots $s$. Here
the sign $\epsilon (s)=(-1)^n$ if $s$ is a sum of
$n$ pairwise perpendicular and perpendicular to $\Lambda$
imaginary simple roots $r_i$, $i \in I$
(here imaginary means $(h_i, h_i)\le 0$),
which are distinct for $i \in I-\tau_0$ ; and $\epsilon (s)=0$ otherwise.
\endproclaim

\proclaim{Statement 6.6$^\prime$}
$$
\cha L(\Lambda )=
\sum_{w\in W}
{\epsilon (w)w(S_\Lambda )} \Big/
e(\rho) \hskip-5pt \prod_{\alpha \in \Delta_+}
{(1-e(-\alpha ))^{\mult~\alpha}}
$$
where $S_\Lambda $ is the sum
$$
S_\Lambda = e(\Lambda+\rho)\sum{\epsilon (s)e(-\pi(s))}
$$
over all sums of imaginary simple roots $s$. Here
$\epsilon (s)=(-1)^{n_\0o}$ if $s$ is a sum of
$n_\0o+n_\1o$
pairwise perpendicular and perpendicular to $\Lambda$
$n_\0o$ imaginary simple roots $r_i$, $i \in I-\tau$, and $n_\1o$
imaginary simple roots $r_i$, $i \in \tau$, which are distinct for
$i \in I-\tau_0$; and $\epsilon (s)=0$ otherwise.
\endproclaim

The corresponding denominator identities are

\proclaim{Statement 6.7$^\prime$}
$$
e(\rho )\prod_{\alpha \in \Delta_+}
{(1-e(-\alpha))^{\mult_\0o \alpha}}
\prod_{\alpha \in \Delta_+}
{(1+e(-\alpha))^{\mult_\1o \alpha }}=
\sum_{w \in W}{\epsilon (w)w(S)}
$$
where $S$ is the sum
$$
S = e(\rho)\sum{\epsilon (s)e(-\pi(s))}
$$
over all sums of imaginary simple roots $s$. Here
the sign $\epsilon (s)=(-1)^n$ if $s$ is a sum of
$n$ pairwise perpendicular
imaginary simple roots $r_i$, $i \in I$,
which are distinct for $i \in I-\tau_0$; and $\epsilon (s)=0$ otherwise.
\endproclaim

\proclaim{Statement 6.8$^\prime$}
$$
e(\rho) \prod_{\alpha \in \Delta_+}
{(1-e(-\alpha ))^{\mult~\alpha}}
= \sum_{w\in W}{\epsilon (w)w(S)}
$$
where $S$ is the sum
$$
S = e(\rho)\sum{\epsilon (s)e(-\pi(s))}
$$
over all sums of imaginary simple roots $s$. Here
$\epsilon (s)=(-1)^{n_\0o}$ if $s$ is a sum of
$n_\0o+n_\1o$ pairwise perpendicular
$n_\0o$ imaginary simple roots $r_i$,
$i \in I-\tau$, and $n_\1o$ imaginary simple roots $r_i$,
$i \in \tau$, which are distinct for
$i \in I-\tau_0$; and $\epsilon (s)=0$ otherwise.
\endproclaim

\vskip5pt

We use the last statement in \S~3 and \S~5. Notation here and in \S~3,
\S~5 are related as follows: $H=M_{II}\otimes \br$,
$I={}_s\Delta$, $\tau ={}_s\Delta^{\im}_{\1o}$,
$G=\geg^{\prime\prime}(M_{II}, {}_s\Delta )$.

\enddocument

\end